# HETEROCLINIC SOLUTIONS FOR A GENERALIZED FRENKEL-KONTOROVA MODEL BY MINIMIZATION METHODS OF RABINOWITZ AND STREDULINSKY

WEN-LONG LI* AND XIAOJUN CUI

ABSTRACT. We study heteroclinic solutions of a generalized Frenkel-Kontorova model. Using the methods of Rabinowitz and Stredulinsky, we prove that if the rotation vector of the configuration is rational and if there is an adjacent pair of periodic configurations, then there is a solution that is heteroclinic in one fixed direction and periodic in other directions. Furthermore, if the above heteroclinic solutions have an adjacent pair, then there is a solution that is heteroclinic in two directions and periodic in other directions. The procedure can be repeated to produce more complex solutions. Thus we obtain a variational construction for these minimal and Birkhoff solutions.

## 1. INTRODUCTION

In this paper, we study a generalized Frenkel-Kontorova (FK, for short) model. To introduce our results, we review some related theories.

**1.1. 1-dimensional FK model and Aubry-Mather theory.** In its simplest form, the 1-dimensional (1-D, for short) FK model describes the motion of a chain of interacting particles ("atoms") subjected to an external on-site periodic potential. A typical example of 1-D FK model can be expressed in the following equation (cf. [12, 7, 8] etc.):

$$\text{(1.1)} \qquad \frac{\mathrm{d}^2 u}{\mathrm{d}t^2}(i) - [u(i+1) + u(i-1) - 2u(i)] + V'(u(i)) = 0,$$

for all $i \in \mathbb{Z}$, where $u(i) \in \mathbb{R}$ is the position of the $i$th particle in the chain. Here $u$ represents the states of the particles of the chain. We call $u : \mathbb{Z} \to \mathbb{R}$ a (lattice) configuration. The external potential $V = V(u) \in C^2(\mathbb{R}, \mathbb{R})$ is a 1-periodic function. (1.1) is difficult to solve since it is not a local problem and it consists of infinitely many equations. In physics, equilibrium states of FK model are of particular concern. A configuration $u$ is said to be a equilibrium state of FK model if it satisfies

$$\text{(1.2)} \qquad -[u(i+1) + u(i-1) - 2u(i)] + V'(u(i)) = 0,$$

for all $i \in \mathbb{Z}$. We refer to [7, 8] for more background and applications of $1-$D FK model.

Among other developments, a breakthrough in the study of FK model is [1, 2]. Almost the same time, J. Mather [13] obtained similar results in monotone twist maps of the annulus. Now their results are named by Aubry-Mather theory. This theory gives a classification on the minimal configurations of a class of Hamiltonians (see [5] for general hypotheses on the Hamiltonians). For (1.2), the Hamiltonian is

$$H(u) = \sum_{i \in \mathbb{Z}} \big\{ \frac{1}{2}[u(i+1) - u(i)]^2 + V(u(i)) \big\}.$$

A configuration is called minimal if for any $v \in \mathbb{R}^{\mathbb{Z}}$ with $v \neq u$ on a finite set, we have $H(u) \leq H(v)$. Obviously, minimal configurations are solutions of (1.2). An important feature of minimal configuration is that it has a rotation number. If

(1.3) $$\lim_{|i| \to \infty} \frac{u(i)}{i} \text{ exists,}$$

denoting this limit by $\alpha$, we say $u$ has rotation number $\alpha$. It can be proved that minimal configuration must have a rotation number (cf. e.g., [5]).

For minimal configurations of FK model (1.2), we briefly introduce Aubry-Mather theory as follows. For any $\alpha \in \mathbb{R}$, the set of minimal configurations with rotation number $\alpha$, denoted by $\mathcal{M}_\alpha$, is not empty. If $\alpha \in \mathbb{R} \setminus \mathbb{Q}$, $\mathcal{M}_\alpha$ is an ordered set and it contains a minimal recurrent set $\mathcal{M}_\alpha^{rec}$. $\mathcal{M}_\alpha^{rec}$ either is $\mathbb{R}$ or is a Cantor set. If $\alpha \in \mathbb{Q}$, $\mathcal{M}_\alpha$ is not an ordered set and it consists of periodic configurations and heteroclinic configurations. The set of periodic configurations of $\mathcal{M}_\alpha$ is ordered. Suppose $u, v \in \mathcal{M}_\alpha$. If there does not exist other periodic configuration $w$ such that $u \leq w \leq v$, then there is a heteroclinic configuration $w_1$ (resp. $w_2$) satisfying $|w_1(i) - u(i)| \to 0$ (resp. $|w_2(i) - v(i)| \to 0$) as $i \to -\infty$, and $|w_1(i) - v(i)| \to 0$ (resp. $|w_2(i) - u(i)| \to 0$) as $i \to \infty$. The set of configurations of $\mathcal{M}_\alpha$ like $w_1$ (resp. $w_2$) is denoted by $\mathcal{M}_{\alpha+}$ (resp. $\mathcal{M}_{\alpha-}$). Aubry-Mather theory tells us that Periodic configurations and heteroclinic configurations in $\mathcal{M}_{\alpha+}$ (resp. $\mathcal{M}_{\alpha-}$) make up an ordered set. For a good survey on this topic, we refer to [5].

After the establishment of Aubry-Mather theory, there are many attempts to generalize it to higher dimensions. One of them is Moser-Bangert theory.

**1.2. Moser-Bangert theory.** In 1986, J. Moser [15] began to generalize Aubry-Mather theory to the case of codimension 1. He considered a nonlinear variational problem on a torus. Under some elliptic conditions, Moser proved that there are minimal solutions of this variational problem. For this variational problem, a function $u \in W^{1,2}_{loc}(\mathbb{R}^n)$ is said to be minimal if it is perturbed by any compact support function, the energy (or Lagrangian, or functional) will not decrease. Of course, if $u$ is minimal, $u$ is a solution of this variational problem. To establish the results similar to Aubry-Mather theory, Moser considered another condition, i.e., without self-intersections (WSI, for short). For any $\mathbf{j} \in \mathbb{Z}^n$ and for any $k \in \mathbb{Z}$, if $u$ satisfies $u(x + \mathbf{j}) - u(x) - k$ does not change sign, then $u$ is WSI. Moser proved that under some elliptic conditions, for every minimal and WSI solution there is a rotation vector $\alpha \in \mathbb{R}^n$, such that $|u(x) - \alpha \cdot x|$ is bounded on $\mathbb{R}^n$. He called $\alpha$ the rotation vector of $u$. Moreover, for any $\alpha \in \mathbb{R}^n$, Moser proved that the set of minimal and WSI solutions with rotation vector $\alpha$, denoted by $\mathcal{M}_\alpha$, is not empty.

V. Bangert [6, 4, 3] made further developments on Moser's problem (known as Moser-Bangert theory). He ([4]) proved that if $\alpha$ is rationally independent, $\mathcal{M}_\alpha$ is an ordered set. The graphs of functions in $\mathcal{M}_\alpha^{rec}$, the minimal recurrent set of $\mathcal{M}_\alpha$, constitute a foliation or lamination. If $\alpha$ is not rationally independent, Bangert ([6]) introduced secondary invariants to classify $\mathcal{M}_\alpha$. Roughly speaking, at this case, $\mathcal{M}_\alpha$ can be decomposed into some ordered sets. Each of these ordered sets is laminated or foliated by periodic solutions and heteroclinic solutions that correspond to secondary invariants.

In a series papers [17, 18, 19, 20, 21] and in the book [22], P. H. Rabinowitz and E. W. Stredulinsky studied an Allen-Cahn type equation which belonged to the variational problem of Moser and Bangert. They used pure variational methods obtaining lots of heteroclinic and homoclinic solutions of the Allen-Cahn equation. These new solutions are not minimal and WSI, but are local minimal. The new results of Rabinowitz and Stredulinsky are based on a new viewpoint on Bangert's heteroclinic solutions. These results are also part of Moser-Bangert theory ([22]).

The relation between Moser-Bangert theory and $n$-dimensional ($n$-D, for short) FK model are explained in the next subsection.

**1.3. $n$-D FK model.** A natural extension of 1-D FK model is $n$-D FK model with $n \geq 2$. Similar to (1.2), $n$-D FK model can be described by the following equation:

$$-(\Delta_D u)(\mathbf{i}) + V'(u(\mathbf{i})) = 0 \tag{1.4}$$

for all $\mathbf{i} \in \mathbb{Z}^n$. Here $V \in C^2(\mathbb{R}, \mathbb{R})$ is 1-periodic and $\Delta_D : \mathbb{R}^{\mathbb{Z}^n} \to \mathbb{R}^{\mathbb{Z}^n}$ is defined as

$$(\Delta_D u)(\mathbf{i}) = \frac{1}{2n} \sum_{\mathbf{j}:\|\mathbf{j}-\mathbf{i}\|=1} (u(\mathbf{j}) - u(\mathbf{i})),$$

where $\|\mathbf{i}\| := \sum_{k=1}^n |\mathbf{i}_k|$. Similar to (1.1), (1.4) is the equilibrium equation of the following equation:

$$\frac{d^2 u}{dt^2}(\mathbf{i}) - (\Delta_D u)(\mathbf{i}) + V'(u(\mathbf{i})) = 0$$

for all $\mathbf{i} \in \mathbb{Z}^n$. A function $u$ defined on $\mathbb{Z}^n$ is also called a (lattice) configuration. Throughout this paper, we use $\mathbf{i}, \mathbf{j}, \mathbf{k}, \mathbf{l}$, etc. (resp. $i, j, k, l$ etc.) to denote elements in $\mathbb{Z}^n$ (resp. $\mathbb{Z}$). Denote by $\mathbf{e}_j$ the vector $(0, \cdots, 1, \cdots, 0)$, i.e., the $j$th component is 1 and the others are 0. If we set $S_{\mathbf{j}}(u) = V(u(\mathbf{j})) + \frac{1}{8n} \sum_{\mathbf{k}:\|\mathbf{k}-\mathbf{j}\|=1} (u(\mathbf{k}) - u(\mathbf{j}))^2$, then (1.4) is the Euler-Lagrange equation

$$\sum_{\mathbf{j} \in \mathbb{Z}^n} \partial_{\mathbf{i}} S_{\mathbf{j}}(u) = \sum_{\mathbf{j}:\|\mathbf{j}-\mathbf{i}\|\leq 1} \partial_{\mathbf{i}} S_{\mathbf{j}}(u) = 0 \text{ for all } \mathbf{i} \in \mathbb{Z}^n \tag{1.5}$$

of the formal sum

$$W(u) := \sum_{\mathbf{j} \in \mathbb{Z}^n} S_{\mathbf{j}}(u). \tag{1.6}$$

Note that (1.6) is a formal sum since the sum may be not convergent. But (1.5) always make sense because the sum has only finite terms.

(1.4) was considered in [9, 24] (in some general forms). Since the dimension $n \geq 2$, the authors also used the property of WSI. The closed relation between Moser-Bangert theory and FK model was explained further in [23].

Following [16, 14] (see also [10, 11] for more general models), we consider a generalized FK model which will be stated explicitly in Section 2. Roughly speaking, we study (1.5) with $S_{\mathbf{j}}$ satisfying some conditions. To attack (1.5), it is natural to use Moser-Bangert theory ([16, 14]). In [16], periodic solutions of (1.5) are obtained. The case for rationally independent rotation vector has also been studied in [16] and the results is analogous to Aubry-Mather theory and Moser-Bangert theorey. In [14], using Bangert's method, the authors defined the secondary invariants and gave a classification of minimal and Birkhoff solutions corresponding to rationally dependent rotation vector. Owing to the variational structure, we can use pure variational method of Rabinowitz and Stredulinsky to obtain similar results of [14] and expect more complex solutions. Now the main results of this paper can be stated.

**1.4. Main results of this paper.** The main results of this paper are as follows. Suppose $\alpha \in \mathbb{Q}^n$. For the generalized FK model of Section 2, we have

- there are periodic solutions with rotation vector $\alpha$;
- if there are an adjacent pair in periodic solutions (gap condition), there exists heteroclinic solutions lying between the adjacent pair, such that they are heteroclinic in one direction and periodic in others;
- solutions heteroclinic in more directions can be obtained provided more gaps conditions.

Comparing to the Allen-Cahn equation considered in [22], the problem in our setting is not local. More care should be taken in applying Rabinowitz and Stredulinsky's methods.

This paper serves in two purposes. On the one hand, it gives a pure variational viewpoint of the heteroclinic solutions of [14]. On the other hand, more homoclinic and heteroclinic solutions are expected to be obtained by pure variational method using heteroclinic solutions of this present paper as building blocks.

This paper is organized as follows. Section 2 gives some preliminaries. In Section 3, we construct heteroclinic solutions that is heteroclinic in $\mathbf{i}_1$ and periodic in $\mathbf{i}_2, \cdots, \mathbf{i}_n$. Then solutions heteroclinic in $\mathbf{i}_1, \mathbf{i}_2$ and periodic in $\mathbf{i}_3, \cdots, \mathbf{i}_n$ are obtained in Section 4. Section 5 includes there generalizations.

## 2. Preliminary

We introduce our generalized FK model. Fix $r \in \mathbb{N}$ and let $B_\mathbf{0}^r = \{\mathbf{k} \in \mathbb{Z}^n \mid \|\mathbf{k}\| \leq r\}$. Assume that $s \in C^2(\mathbb{R}^{B_\mathbf{0}^r}, \mathbb{R})$ satisfies

(S1) $s(u + 1_{B_\mathbf{0}^r}) = s(u)$, where $1_{B_\mathbf{0}^r}$ is the constant function 1 on $B_\mathbf{0}^r$;

(S2) $s$ is bounded from below and coercive in the following sence,

$$\lim_{|u(\mathbf{k}) - u(\mathbf{j})| \to \infty} s(u) = \infty, \text{ for } \mathbf{k}, \mathbf{j} \in B_\mathbf{0}^r \text{ with } \|\mathbf{k} - \mathbf{j}\| = 1;$$

(S3) $\partial_{\mathbf{k},\mathbf{j}} s \leq 0$ for $\mathbf{k}, \mathbf{j} \in B_\mathbf{0}^r$ with $\mathbf{k} \neq \mathbf{j}$, while $\partial_{\mathbf{0},\mathbf{j}} s < 0$ for $\|\mathbf{j}\| = 1$.

For $u \in \mathbb{R}^{\mathbb{Z}^n}$, set $S_\mathbf{j}(u) = s(\tau^n_{-\mathbf{j}_n} \cdots \tau^1_{-\mathbf{j}_1} u|_{B_\mathbf{0}^r})$, where $\tau^j_{-k} : \mathbb{R}^{\mathbb{Z}^n} \to \mathbb{R}^{\mathbb{Z}^n}$ is defined by $\tau^j_{-k} u(\mathbf{i}) = u(\mathbf{i} + k\mathbf{e}_j)$. With these local potential $S_\mathbf{j}$, we can define the formal sum (1.6) and its Euler-Lagrange equation

$$(2.1) \quad \sum_{\mathbf{j} \in \mathbb{Z}^n} \partial_\mathbf{i} S_\mathbf{j}(u) = \sum_{\mathbf{j}: \|\mathbf{j} - \mathbf{i}\| \leq r} \partial_\mathbf{i} S_\mathbf{j}(u) = 0 \text{ for all } \mathbf{i} \in \mathbb{Z}^n.$$

Sometimes it is useful to consider $u \in \mathbb{R}^{B_\mathbf{0}^r}$ as a vector $(u(\mathbf{i}))_{\mathbf{i} \in B_\mathbf{0}^r} \in \mathbb{R}^{\#B_\mathbf{0}^r}$. Here $\#B_\mathbf{0}^r$ denote the cardinality of $B_\mathbf{0}^r$. Thus $s \in C^2(\mathbb{R}^{B_\mathbf{0}^r}, \mathbb{R})$ can be seen as a function of several variables. So $\|s\|_{L^\infty(\mathbb{R}^{\#B_\mathbf{0}^r})}$ is well-defined. These viewpoints will be useful in our analysis. See Propositions 3.1, 3.16, Remarks 3.17, 3.19, etc.

Before going further we recall some definitions. A solution $u$ of (2.1) is a configuration defined on $\mathbb{Z}^n$ satisfying (2.1). It is equivalent to the stationary point of the local potentials $S_\mathbf{j}$ (for the definition of stationary point, cf. [16, Definition 2.2]). On the lattice $\mathbb{Z}^n$ we define interior points of a subset as follows. For any subset $B$ of $\mathbb{Z}^n$, the $r$-interior of $B$, denoted by $int_r(B)$, is defined by $int_r(B) = \{\mathbf{i} \in B \mid B_\mathbf{i}^r \subset B\}$, where $B_\mathbf{i}^r := \{\mathbf{k} \in \mathbb{Z}^n \mid \|\mathbf{k} - \mathbf{i}\| \leq r\}$ (cf. [14, p.1525, line 15]). In this paper, the main object is minimal and Birkhoff configurations whose definitions we now introduce.

**Definition 2.1** (cf. [16, Definition 2.3]). *A configuration $u : \mathbb{Z}^n \to \mathbb{R}$ is called a (global) minimizer for the potentials $S_\mathbf{j}$ (or for potential $s$) if for every finite subset $B \subset \mathbb{Z}^n$ and every $v : \mathbb{Z}^n \to \mathbb{R}$ with support in $int_r(B)$,*

$$W_B(u+v) - W_B(u) = \sum_{\mathbf{j} \in B}(S_\mathbf{j}(u+v) - S_\mathbf{j}(u)) \geq 0,$$

*where the support of $v$ is $supp(v) := \{\mathbf{i} \in \mathbb{Z}^n \,|\, v(\mathbf{i}) \neq 0\}$ and $W_B : \mathbb{R}^{\mathbb{Z}^n} \to \mathbb{R}$ is defined as*

$$W_B(u) = \sum_{\mathbf{j} \in B} S_\mathbf{j}(u).$$

To define the Birkhoff configuration, we introduce three partial order relations on $\mathbb{R}^{\mathbb{Z}^n}$.

**Definition 2.2** (cf. [16, Definition 3.2]). *We define the relations $\leq, \lneqq$ and $<$ on $\mathbb{R}^{\mathbb{Z}^n}$ by:*
- *$u \leq v$ if $u(\mathbf{i}) \leq v(\mathbf{i})$ for every $\mathbf{i} \in \mathbb{Z}^n$;*
- *$u \lneqq v$ if $u \leq v$ and there is some $\mathbf{i} \in \mathbb{Z}^n$ such that $u(\mathbf{i}) \neq v(\mathbf{i})$;*
- *$u < v$ if $u(\mathbf{i}) < v(\mathbf{i})$ for every $\mathbf{i} \in \mathbb{Z}^n$.*

*The relations $\geq, \gneqq$ and $>$ are defined similarly.*

**Definition 2.3** (cf. [14, Definition 2.1], [22, p.3, line 25]). *A configuration $u$ is said to be Birkhoff if $\{\tau_j^k u \,|\, j \in \mathbb{Z} \text{ and } 1 \leq k \leq n\}$ is totally ordered, i.e., for all $j \in \mathbb{Z}$ and $1 \leq k \leq n$, it follows that*

$$\tau_j^k u < u, \quad or \quad \tau_j^k u = u, \quad or \quad \tau_j^k u > u.$$

*If $u \leq \tau_{-j}^k u$ (resp. $u < \tau_{-j}^k u$) holds, we say $u$ is $j$-monotone (resp. strictly $j$-monotone) in $\mathbf{i}_k$.*

As in Aubry-Mather theory, an important feature of Birkhoff configuration is that it has a rotation vector. Note in Aubry-Mather theory, i.e., in the $1-D$ case, the property of minimal implies Birkhoff, cf. [5].

**Definition 2.4** (cf. [16, Definition 3.1]). *Let $u : \mathbb{Z}^n \to \mathbb{R}$. $\alpha \in \mathbb{R}^n$ is said to be the rotation vector of $u$ if for all $\mathbf{i} \in \mathbb{Z}^n$, the limit*

$$\lim_{|m| \to \infty} \frac{u(m\mathbf{i})}{m} \text{ exists and is equal to } \langle \alpha, \mathbf{i} \rangle.$$

**Lemma 2.5** (cf. [16, Lemma 3.5]). *Let $u : \mathbb{Z}^n \to \mathbb{R}$ be a Birkhoff configuration. Then $u$ has a rotation vector $\alpha = \alpha(u)$ and*

$$|u(\mathbf{i}) - u(\mathbf{0}) - \langle \alpha(u), \mathbf{i} \rangle| \leq 1.$$

The ordered relations among solutions of (2.1) are the key point of our analysis. We have:

**Lemma 2.6** (cf. [14, Lemma 2.5]; [16, Lemma 4.5]). *Assume that $u$ and $v$ are solutions of (2.1) and $u \leq v$. Then $u < v$ or $u = v$.*

The deduction of the next corollary appears repeatedly, so we pick it out as a corollary for convenience. The idea of the proof follows from, for example, [22, Proposition 2.2].

**Corollary 2.7.** *Assume that $u, v$ are solutions of (2.1). If $\psi := \min(u, v)$ or $\phi := \max(u, v)$ is a solution of (2.1), then*

$$u < v, \quad or \quad u = v, \quad or \quad u > v.$$

**Proof:** Suppose $\psi$ is a solution of (2.1). Since $\psi \leq u$, by Lemma 2.6, $\psi < u$ or $\psi = u$. If $\psi < u$, then $v = \psi < u$. If $\psi = u$, then $u = \psi \leq v$. Using Lemma 2.6 again yields $u < v$ or $u = v$. The case for $\phi$ can be proved similarly. □

**Lemma 2.8** (cf. [14, Lemma 2.6]). *For $u, v \in \mathbb{R}^{\mathbb{Z}^n}$ and an arbitrary finite set $B \subset \mathbb{Z}^n$, we have*
$$W_B(u) + W_B(v) \geq W_B(\phi) + W_B(\psi),$$
*where $\phi, \psi$ are defined by $\phi = \max(u,v)$, $\psi = \min(u,v)$.*

A variant of Lemma 2.8 is often used conjunctively with Corollary 2.7 to obtain the order relation of $u, v$. For the proof of Lemma 2.8, we refer the reader to [14]. The next is the convergence we need in this paper.

**Definition 2.9.** *Let $E \subset \mathbb{Z}^n$. By saying $u_k \to u$ in $\mathbb{R}^E$ as $k \to \infty$, we mean $u_k(\mathbf{i}) \to u(\mathbf{i})$ for all $\mathbf{i} \in E$, i.e., the convergence is pointwise. When $u_k \to u$ in $\mathbb{R}^{\mathbb{Z}^n}$, we say $u_k \to u$ pointwisely. By saying $u_k \to u$ in $\|\cdot\|_E$ as $k \to \infty$, we mean*
$$\|u_k - u\|_E := \sum_{\mathbf{j} \in E} |u_k(\mathbf{j}) - u(\mathbf{j})| \to 0.$$

*Of course, for any bounded set $E \subset \mathbb{Z}^n$, $u_k \to u$ in $\mathbb{R}^E$ is equivalent to $u_k \to u$ in $\|\cdot\|_E$.*

**Remark 2.10.** *Definition 2.9 provides an approach to use the method of Rabinowitz and Stredulinsky. The norm $\|\cdot\|_E$ will replace $\|\cdot\|_{L^2(E)}$ and thus $\|\cdot\|_{W^{1,2}(E)}$ since we do not have the term $\|\nabla u\|_{L^2(E)}$ in our setting.*

A configuration $u$ is called to have bounded action if there exists $C > 0$, such that $|u(\mathbf{k}) - u(\mathbf{j})| \leq C$ for all $\mathbf{k}, \mathbf{j} \in \mathbb{Z}^n$ with $\|\mathbf{k} - \mathbf{j}\| = 1$ (cf. [14, p.1525, line -3]). The following lemma provides an important estimate, which will be used on the functionals $J_1, J_2$, etc. (See Sections 3, 4 for the definitions.)

**Lemma 2.11** (cf. [14, Lemma 2.4]). *Assume that $u$ and $v$ have bounded action with bounded constant $C$. Then there exists a constant $L = L(C,r) > 0$ such that for each finite set $B \subset \mathbb{Z}^n$,*
$$|W_B(u) - W_B(v)| \leq L \sum_{\mathbf{i} \in \bar{B}} |u(\mathbf{i}) - v(\mathbf{i})|,$$
*where $\bar{B} = \cup_{\mathbf{j} \in B} B_{\mathbf{j}}^r = \cup_{\mathbf{j} \in B} \{\mathbf{k} \in \mathbb{Z}^n \mid \|\mathbf{k} - \mathbf{j}\| \leq r\}$.*

In [16] the authors had constructed Aubry-Mather sets associated to every rotational vector $\alpha \in \mathbb{R}^n$. Miao, et al. ([14]) used Bangert's idea ([6]) proving the existence of configurations corresponding to the secondary invariants. We want to prove similar results as that in [14] by Rabinowitz-Stredulinsky method, that is, minimization method. We begin with constructing Aubry-Mather set corresponding to $\alpha = \mathbf{0}$.

We say a function $u : \mathbb{Z} \to \mathbb{R}$ is of period 1 if $u(i+1) = u(i)$ for all $i \in \mathbb{Z}$. Denote the set of all 1-periodic functions defined on $\mathbb{Z}$ by $\mathbb{R}^{\mathbb{Z}/\{1\}}$. Similarly, we can define the set $\mathbb{R}^{(\mathbb{Z}/\{1\})^n}$, $\mathbb{R}^{\mathbb{Z} \times (\mathbb{Z}/\{1\})^{n-1}}$, etc. Let $J_{\mathbf{0}}(u) = S_{\mathbf{0}}(u)$ and define
$$\Gamma_0 = \mathbb{R}^{(\mathbb{Z}/\{1\})^n},$$

(2.2) $$c_0 = \inf_{u \in \Gamma_0} J_0(u),$$

and
$$\mathcal{M}_0 = \{u \in \Gamma_0 \mid J_\mathbf{0}(u) = c_0\}.$$
It is easily seen that $u \in \Gamma_0$ if and only if $u \equiv t$ for some $t \in \mathbb{R}$ and thus
$$c_0 = \min_{t \in \mathbb{R}} S_\mathbf{0}(t) = \min_{0 \leq t \leq 1} S_\mathbf{0}(t).$$
Hence the constant configuration $t_0 \pm j \in \mathcal{M}_0$ for all $j \in \mathbb{Z}$, where $t_0 \in [0, 1)$ satisfies $S_\mathbf{0}(t_0) = \min_{0 \leq t \leq 1} S_\mathbf{0}(t)$. Moreover, imbedding $\mathcal{M}_0$ to $\mathbb{R}$, we have

**Theorem 2.12.** *$\mathcal{M}_0$ is a nonempty ordered set and the elements in $\mathcal{M}_0$ are solutions of* (2.1).

The first assertion is easy to prove. The second statement follows from [16, Theorem 4.8]. Suppose

($*_0$) there are adjacent $v_0, w_0 \in \mathcal{M}_0$ with $v_0 < w_0$.

Throughout this paper, "there are adjacent $v, w \in A$ with $v < w$" means there does not exist $u \in A$ satisfying $v \lneq u \lneq w$. In Aubry-Mather theory, condition like ($*_0$) is a sufficient condition for the existence of heteroclinic solutions. Rabinowitz and Stredulinsky construct heteroclinic solutions under ($*_0$) and we will adopt their notations. ($*_0$) can be easily fulfilled, for instance, when $S_\mathbf{0}(t)$ has finite minimal points in $[0, 1]$. ($*_0$) is a generic property. See Proposition 3.18 below.

## 3. Solutions heteroclinic in $\mathbf{i}_1$

In this section, under ($*_0$) we establish the solutions heteroclinic from $v_0$ to $w_0$ in $\mathbf{i}_1$. Let $v, w \in \mathcal{M}_0$ with $v < w$. At this moment, we do not require $v, w$ are adjacent in $\mathcal{M}_0$.

Define
$$\hat{\Gamma}_1 = \hat{\Gamma}_1(v, w) = \{u \in \mathbb{R}^{\mathbb{Z} \times (\mathbb{Z}/\{1\})^{n-1}} \mid v \leq u \leq w\}.$$
For $i \in \mathbb{Z}$, set $\mathbf{T}_i = (i, 0, \cdots, 0)$. For $u \in \hat{\Gamma}_1$ and $i \in \mathbb{Z}$, define
$$J_{1,i}(u) = J_\mathbf{0}(\tau^1_{-i} u) - c_0 = S_{\mathbf{T}_i}(u) - c_0$$
with $c_0$ as in (2.2).

For $p, q \in \mathbb{Z}$ with $p \leq q$ and $u \in \hat{\Gamma}_1$, set
$$J_{1;p,q}(u) = \sum_{i=p}^{q} J_{1,i}(u).$$

Let us begin from the study on periodic configurations.

**Proposition 3.1.** *Let $\mathbf{l} \in \mathbb{N}^n$ and*
$$\Gamma_0(\mathbf{l}) = \{u : \mathbb{Z}^n \to \mathbb{R} \mid u(\mathbf{i} + \mathbf{l}_k \mathbf{e}_k) = u(\mathbf{i}), 1 \leq k \leq n\}.$$
*Set*
$$J_0^\mathbf{l}(u) = \sum_{\substack{0 \leq \mathbf{i}_k < \mathbf{l}_k \\ 1 \leq k \leq n}} S_\mathbf{i}(u) \quad \text{and} \quad c_0(\mathbf{l}) = \inf_{u \in \Gamma_0(\mathbf{l})} J_0^\mathbf{l}(u).$$
*Then*
$$\mathcal{M}_0(\mathbf{l}) := \{u \in \Gamma_0(\mathbf{l}) \mid J_0^\mathbf{l}(u) = c_0(\mathbf{l})\} \neq \emptyset.$$
*Moreover, $\mathcal{M}_0(\mathbf{l}) = \mathcal{M}_0$ and $c_0(\mathbf{l}) = (\prod_{i=1}^n \mathbf{l}_i) c_0$.*

**Proof:** It is easy to see that $\mathcal{M}_0(\mathbf{l}) \neq \emptyset$. In fact, by the definition of $\Gamma_0(\mathbf{l})$, $J_0^{\mathbf{l}}$ can be considered as a function of finite variables. Note that if $u \in \Gamma_0(\mathbf{l})$ then so is $u \pm j$ for all $j \in \mathbb{Z}$. We may assume the minimizing sequence $u_n$ satisfying $u_n(\mathbf{0}) \in [0,1]$. Since $J_0^{\mathbf{l}}(u_n)$ is bounded, by (S2) we deduce that $\{u_n(\mathbf{i})\}_n$ is bounded for all $0 \leq \mathbf{i}_k < \mathbf{l}_k$ and $1 \leq k \leq n$. Since

(3.1) $\quad$ for every $a \in \mathbb{R}^{\mathbb{Z}^n}$ the set $\{u \in \mathbb{R}^{\mathbb{Z}^n} \mid |u(\mathbf{i})| \leq a(\mathbf{i})$ for all $\mathbf{i} \in \mathbb{Z}^n\}$ is compact

(cf. [5, (1.1)]), minimization method ensures $\mathcal{M}_0(\mathbf{l}) \neq \emptyset$ and the elements of $\mathcal{M}_0(\mathbf{l})$ are solutions of (2.1) (cf. [16, Theorem 4.8]).

We claim that $\mathcal{M}_0(\mathbf{l})$ is an ordered set. Suppose $v, w \in \mathcal{M}_0(\mathbf{l})$. Set $\phi = \max(v, w)$, $\psi = \min(v, w)$. Then $\phi, \psi \in \Gamma_0(\mathbf{l})$ and by Lemma 2.8,

(3.2) $$J_0^{\mathbf{l}}(\phi) + J_0^{\mathbf{l}}(\psi) \leq J_0^{\mathbf{l}}(v) + J_0^{\mathbf{l}}(w) = 2c_0(\mathbf{l}).$$

Since
$$J_0^{\mathbf{l}}(\phi), J_0^{\mathbf{l}}(\psi) \geq c_0(\mathbf{l}),$$
(3.2) implies $J_0^{\mathbf{l}}(\phi) = J_0^{\mathbf{l}}(\psi) = c_0(\mathbf{l})$, so $\phi, \psi \in \mathcal{M}_0(\mathbf{l})$. Therefore $\phi$ and $\psi$ are solutions of (2.1). By Corollary 2.7, $v < w$ or $v = w$ or $v > w$. Hence $\mathcal{M}_0(\mathbf{l})$ is an ordered set.

Now to prove $\mathcal{M}_0(\mathbf{l}) = \mathcal{M}_0$ and $c_0(\mathbf{l}) = (\prod_{i=1}^n \mathbf{l}_i) c_0$, it suffices to prove:

(3.3) $\quad$ if $u \in \mathcal{M}_0(\mathbf{l})$ then $u(\mathbf{i} + \mathbf{e}_j) = u(\mathbf{i})$, $\quad 1 \leq j \leq n$.

The proof of (3.3) is the same to [22, (2.6)]. Here we give the proof for the sake of completeness. Suppose $u \in \mathcal{M}_0(\mathbf{l})$. Noting $u(\mathbf{i} + \mathbf{e}_j) \in \mathcal{M}_0(\mathbf{l})$, $j = 1, \ldots, n$, and $\mathcal{M}_0(\mathbf{l})$ is ordered, we have either (3.3) holds or

(3.4) $\quad$ (i) $u(\mathbf{i} + \mathbf{e}_j) > u(\mathbf{i})$ $\quad$ or $\quad$ (ii) $u(\mathbf{i} + \mathbf{e}_j) < u(\mathbf{i})$

for each $j$. But if (3.4) (i) is satisfied, we obtain
$$u(\mathbf{i}) = u(\mathbf{i} + \mathbf{l}_j \mathbf{e}_j) \geq \cdots \geq u(\mathbf{i} + \mathbf{e}_j) > u(\mathbf{i}),$$
which is a contradiction. Similarly (3.4) (ii) will not hold and then (3.3) is proved. The proof of Proposition 3.1 is complete. $\quad\square$

With Proposition 3.1 in hands, we obtain a lower bound for $J_{1;p,q}(u)$.

**Proposition 3.2.** *If $u \in \hat{\Gamma}$ and $p < q \in \mathbb{Z}$, there is a constant $K_1 = K_1(v,w) \geq 0$, such that*
$$J_{1;p,q}(u) \geq -K_1.$$

**Proof:** First we add an additional condition $q - p \geq 2r + 2$. Taking $u \in \hat{\Gamma}_1$, define

(3.5) $$\chi = \begin{cases} u, & p+1 \leq \mathbf{i}_1 \leq q-1, \\ v, & \mathbf{i}_1 = p, q, \end{cases}$$

and extend $\chi$ as a $(q + 1 - p)$-periodic function of $\mathbf{i}_1$. By Proposition 3.1,
$$0 \leq J_{1;p,q}(\chi) = J_{1;p,p+r}(\chi) + J_{1;p+r+1,q-r-1}(u) + J_{1;q-r,q}(\chi),$$
or

(3.6) $$J_{1;p+r+1,q-r-1}(u) \geq -J_{1;p,p+r}(\chi) - J_{1;q-r,q}(\chi).$$

Since $u \in \hat{\Gamma}_1$, by Lemma 2.11 with $C = w(\mathbf{0}) - v(\mathbf{0})$, there is an $L = L(C,r) > 0$ such that for any $i \in \mathbb{Z}$,
$$|J_{1,i}|(\chi) \leq \sum_{\mathbf{j} \in \bar{\mathbf{T}}_i} L(\chi - v)(\mathbf{j}) \leq L(\#B_0^r)(w-v)(\mathbf{0}),$$

where $\#A$ is the cardinality of the set $A$. Thus (3.6) implies

(3.7) $$J_{1;p+r+1,q-r-1}(u) \geq -(2r+2)L(\#B_{\mathbf{0}}^r)(w-v)(\mathbf{0}).$$

Note that $s \geq -M$ for some $M > 0$ by (S2), then

$$J_{1;p,p+r}(u) \geq (r+1)(-M-|c_0|)$$

and

$$J_{1;q-r,q}(u) \geq (r+1)(-M-|c_0|).$$

Hence

$$\begin{aligned}&J_{1;p,q}(u)\\ =&J_{1;p,p+r}(u)+J_{1;p+r+1,q-r-1}(u)+J_{1;q-r,q}(u)\\ \geq& -(r+1)(M+|c_0|)-(2r+2)L(\#B_{\mathbf{0}}^r)(w-v)(\mathbf{0})-(r+1)(M+|c_0|)\\ =& -(2r+2)[|c_0|+M+L(\#B_{\mathbf{0}}^r)(w-v)(\mathbf{0})].\end{aligned}$$

Thus we complete the proof of Proposition 3.2 with the additional condition that

$$q-p \geq 2r+2.$$

Now for $0 \leq q-p \leq 2r+1$, $J_{1;p,q}(u) \geq -(2r+1)(M+|c_0|)$. Letting

$$K_1 = (2r+2)\bigl[|c_0|+M+L(\#B_{\mathbf{0}}^r)(w-v)(\mathbf{0})\bigr],$$

we complete the proof of Proposition 3.2. □

Now following [22] we define $J_1(u)$ as

$$J_1(u) = \liminf_{\substack{p\to-\infty\\q\to\infty}} J_{1;p,q}(u)$$

for $u \in \hat{\Gamma}_1$. An upper bound for $J_{1;p,q}(u)$ is ready (cf. [22, Lemma 2.22]).

**Lemma 3.3.** *If $u \in \hat{\Gamma}_1$ and $p \leq q \in \mathbb{Z}$, then*

$$J_{1;p,q}(u) \leq J_1(u) + 2K_1.$$

The proof of Lemma 3.3 is exactly same to that of [22, Lemma 2.22], so we omit it.

For $E \subset \mathbb{Z}^n$ and $u \in \mathbb{R}^{\mathbb{Z}^n}$, let $\|u\|_E = \sum_{\mathbf{i} \in E} |u(\mathbf{i})|$. Define

$$\Gamma_1 = \Gamma_1(v,w) = \{u \in \hat{\Gamma}_1 \mid \|u-v\|_{\mathbf{T}_i} \to 0, i \to -\infty; \|u-w\|_{\mathbf{T}_i} \to 0, i \to \infty\}.$$

**Proposition 3.4.** *If $u \in \Gamma_1$, then*

(3.8) $$\|\tau^1_{-i}u - v\|_{\mathbf{T}_0} \to 0, \quad i \to -\infty,$$

(3.9) $$\|\tau^1_{-i}u - w\|_{\mathbf{T}_0} \to 0, \quad i \to \infty,$$

(3.10) $$J_{1,i}(u) \to 0, \quad |i| \to \infty.$$

*If $u \in \Gamma_1$ and $J_1(u) < \infty$, then*

(3.11) $$J_1(u) = \lim_{\substack{p\to-\infty\\q\to\infty}} J_{1;p,q}(u).$$

**Proof:** (3.8)-(3.9) easily follow from the definitions of $\Gamma_1$ and $\tau^1_{-i}u$. By (3.8), the periodicity of $u, v$, and the continuity of $S_{\mathbf{0}}$, we have $\lim_{i\to-\infty} J_{1,i}(u) = 0$. Similarly we can prove $\lim_{i\to\infty} J_{1,i}(u) = 0$, so (3.10) follows. To prove (3.11), it suffices to show that

(3.12) $$(i) \lim_{p\to-\infty} J_{1;p,0}(u) \quad \text{and} \quad (ii) \lim_{q\to\infty} J_{1;0,q}(u)$$

exist.

Since the proofs of (i) and (ii) of (3.12) are the same, we will only verify (3.12) (i). Set
$$\mathcal{P} = \{p \in \mathbb{Z} \,|\, p < 0 \text{ and } J_{1,p}(u) \leq 0\}.$$
If the cardinality of $\mathcal{P}$ is finite, i.e., $\#\mathcal{P} < \infty$, $J_{1;p,0}(u)$ is monotone nondecreasing as $p \to -\infty$. Since $J_{1;p,0}(u) \leq J_1(u) + 2K_1$, the proof of (3.11) (i) is complete. Now assume $\#\mathcal{P} = \infty$. Suppose (3.11) (i) is false, i.e., $J_{1;p,0}(u)$ dose not converge as $p \to -\infty$. Let
$$l^- = \liminf_{p \to -\infty} J_{1;p,0}(u), \quad l^+ = \limsup_{p \to -\infty} J_{1;p,0}(u),$$
then $l^+ > l^- \geq -K_1$. Choose $\epsilon$ such that

(3.13) $$(l^+ - l^-)/5 > \epsilon > 0.$$

Now we need a technical lemma.

**Lemma 3.5.** *For any $\gamma > 0$, there is a $\delta = \delta(\gamma) > 0$ such that if $u \in \Gamma_1(v, w)$, $p, q \in \mathbb{Z}$, with $q - p \geq 4r + 2$ and*

(3.14) $$\sum_{i=-r}^{2r} \|u - v\|_{\mathbf{T}_{j+i}} \leq \delta \quad \text{or} \quad \sum_{i=-r}^{2r} \|u - w\|_{\mathbf{T}_{j+i}} \leq \delta$$

*for $j = p$ and $q - r$, then*

(3.15) $$J_{1;p+r+1,q-r-1}(u) \geq -\gamma.$$

*Proof of Lemma 3.5.* Suppose $\sum_{i=-r}^{2r} \|u - v\|_{\mathbf{T}_{j+i}} \leq \delta$ for $j = p$ and $q - r$. Taking $\chi$ as in (3.5), if

(3.16) $$\sum_{i=0}^{r} |J_{1,p+i}(\chi)| + \sum_{i=0}^{r} |J_{1,q-r+i}(\chi)| \leq \gamma,$$

by (3.6) we obtain (3.15). But (3.16) follows from (3.14), the definition of $\chi$, and the continuity of $J_{1,i}$ (in $\mathbb{R}^{\mathbb{Z}^n}$) for $i \in \mathbb{Z}$. The case of $\sum_{i=-2r}^{r} \|u - w\|_{\mathbf{T}_{j+i}} \leq \delta$ can be proved similarly. This complete the proof of Lemma 3.5

Let us continue the proof of Proposition 3.4. Choose $\gamma = \epsilon$ and $\delta = \delta(\epsilon)$ in Lemma 3.5. By (3.10) and (3.8), there is a $p_0 \in \mathcal{P}$ such that for $p \leq p_0$

(3.17) $$\begin{cases} J_{1,p}(u) \geq -\epsilon/(r+1), \\ \|\tau_{-p}^1 u - v\|_{\mathbf{0}} \leq \delta/(3r+1). \end{cases}$$

Thus by Lemma 3.5,

(3.18) $$J_{1;p+r+1,q-r-1}(u) \geq -\epsilon,$$

for $p, q \in \mathbb{Z}$ and $p + 2r + 1 < q - 2r - 1 < q + r + 1 \leq p_0$. Choose two sequences $(p_k), (q_k) \subset -\mathbb{N}$ such that $q_{k+1} + 9r + 4 < p_k + 5r + 2 < q_k + r < p_0$ and
$$J_{1;p_k,0}(u) \to l^-; \quad J_{1;q_k,0}(u) \to l^+, \quad k \to \infty.$$
Hence there exists a $k_0$ such that for $k \geq k_0$,

(3.19) $$J_{1;p_k,0}(u) \leq l^- + \epsilon; \quad J_{1;q_k,0}(u) \geq l^+ - \epsilon.$$

Let $\hat{q}_k = \max\{q \in \mathcal{P} \,|\, q < q_k\}$ and $\hat{p}_k = \min\{p \in \mathcal{P} \,|\, p \geq p_k\}$. We may assume

(3.20) $$\hat{q}_k - \hat{p}_k \geq 4r + 2.$$

If not, replacing $(p_k)$ and $(q_k)$ by their suitable subsequences we have (3.20). Then
$$J_{1;\hat{q}_k+1,q_k-1}(u) \geq 0, \quad J_{1;p_k,\hat{p}_k-1}(u) \geq 0,$$
and $J_{1;\hat{q}_k+1,q_k-1}(u)$ (resp. $J_{1;p_k,\hat{p}_k-1}(u)$) does not exist if $\hat{q}_k = q_k - 1$ (resp. $\hat{p}_k = p_k$). By (3.19),
(3.21) $$J_{1;\hat{p}_k,0}(u) \leq l^- + \epsilon, \quad J_{1;\hat{q}_k+1,0}(u) \geq l^+ - \epsilon.$$
Consequently, by (3.21) and (3.13),
(3.22) $$J_{1;\hat{p}_k,\hat{q}_k}(u) = J_{1;\hat{p}_k,0}(u) - J_{1;\hat{q}_k+1,0}(u) \leq l^- + \epsilon - (l^+ - \epsilon) < -3\epsilon.$$
However, by (3.18),
$$J_{1;\hat{p}_k+r+1,\hat{q}_k-r-1}(u) \geq -\epsilon.$$
With this inequality and by taking $p = \hat{p}_k, \hat{p}_k + 1, \cdots, \hat{p}_k + r, \hat{q}_k - r, \hat{q}_k - r + 1, \cdots, \hat{q}_k$ in (3.17), we get
$$J_{1;\hat{p}_k,\hat{q}_k}(u) \geq -3\epsilon,$$
contrary to (3.22). Thus $l^- = l^+$, which completes the proof of Proposition 3.4. □

**Corollary 3.6.** *Suppose $u \in \hat{\Gamma}_1(v,w)$, $J_1(u) < \infty$, and $u$ is 1-monotone in $\mathbf{i}_1$. Then $u \in \mathcal{M}_0$ or there are $\phi, \psi \in \mathcal{M}_0$ with $v \leq \phi < \psi \leq w$ such that $u \in \Gamma_1(\phi, \psi)$.*

**Proof:** By (3.1), any sequence $(u_k)$ in $\hat{\Gamma}_1(v,w)$ is precompact. Notice that $J_0$ is continuous with respect to pointwise convergence. Then the proof of Corollary 3.6 follows as in [22, Corollary 2.49]. □

Corollary 3.6 highlights the Birkhoff configurations in $\hat{\Gamma}_1(v,w)$. It implies that Birkhoff configuration either is periodic or heteroclinic to adjacent pair of periodic configurations under the mild condition $J_1(u) < \infty$.

To apply minimization argument, besides continuous of the functional, the compact property of minimizing sequences should be considered. Fortunately, in our setting it is easy to verify (at least for $J_1$).

**Proposition 3.7.** *Let $\mathcal{Y} \subset \hat{\Gamma}_1(v,w)$ and define*
(3.23) $$c(\mathcal{Y}) = \inf_{u \in \mathcal{Y}} J_1(u)$$
*Suppose $(u_k)$ is a minimizing sequence for (3.23), then there is a $U \in \hat{\Gamma}_1$ such that along a subsequence, $u_k \to U$ pointwise. If $c(\mathcal{Y}) < \infty$, then*
(3.24) $$-K_1 \leq J_1(U) \leq c(\mathcal{Y}) + 1 + 2K_1,$$
*with $K_1$ as in Proposition 3.2.*

**Proof:** The first assertion is easily obtained by (3.1). Suppose $c(\mathcal{Y}) < \infty$. We may assume $(u_k)$ satisfying $J_1(u_k) \leq c(\mathcal{Y}) + 1$. By Proposition 3.2, Lemma 3.3, and the continuity of $J_{1;p,q}$, we have
(3.25) $$-K_1 \leq J_{1;p,q}(U) \leq c(\mathcal{Y}) + 1 + 2K_1$$
for any $p \leq q$ and (3.24) follows. □

Proposition 3.7 can not be used directly since the limit point $U$ may not belong to $\mathcal{Y}$. So it cannot imply that $U$ is a solution of (2.1). The following proposition provides a criteria to ensure that minimal point obtained by minimizing sequence of $J_1$ on $\Gamma_1$ is a solution of (2.1). We need a notation. For $\mathbf{T}_i = (i, 0, \cdots, 0)$, define $\delta_{\mathbf{T}_i}(\mathbf{j}) = 1$ if $\mathbf{j}_1 = i$, otherwise $\delta_{\mathbf{T}_i}(\mathbf{j}) = 0$.

**Proposition 3.8.** *Let $\mathcal{Y} \subset \hat{\Gamma}_1(v,w)$. If $c(\mathcal{Y}) < \infty$ and there is a minimizing sequence $(u_k)$ for $c(\mathcal{Y})$ such that for some $i \in \mathbb{Z}$, the function $\delta_{\mathbf{T}_i}$ and some $t_0 > 0$, we have*

$$c(\mathcal{Y}) \leq J_1(u_k + t\delta_{\mathbf{T}_i}) + \epsilon_k \tag{3.26}$$

*for all $|t| \leq t_0$, where $\epsilon_k \to 0$ as $k \to \infty$. Then the limit $U$ of $u_k$ satisfies (2.1) at $\mathbf{T}_i$. Moreover, $U$ satisfies (2.1) at any $\mathbf{j}$ with $\mathbf{j}_1 = (\mathbf{T}_i)_1 = i$.*

**Proof:** Suppose $(u_k)$ is the minimizing sequence for (3.23) satisfying (3.26). Define $\epsilon'_k$ via

$$J_1(u_k) = c(\mathcal{Y}) + \epsilon'_k,$$

so $\epsilon'_k \to 0$ as $k \to \infty$. By (3.26),

$$J_1(u_k) = c(\mathcal{Y}) + \epsilon'_k \leq J_1(u_k + t\delta_{\mathbf{T}_i}) + \epsilon_k + \epsilon'_k,$$

Thus by the definition of $\delta_{\mathbf{T}_i}$,

$$\sum_{j:|j-i|\leq r} J_{1,j}(u_k) \leq \sum_{j:|j-i|\leq r} J_{1,j}(u_k + t\delta_{\mathbf{T}_i}) + \epsilon_k + \epsilon'_k.$$

Letting $k \to \infty$ we have

$$\sum_{j:|j-i|\leq r} J_{1,j}(U) \leq \sum_{j:|j-i|\leq r} J_{1,j}(U + t\delta_{\mathbf{T}_i}),$$

or

$$\sum_{j:|j-i|\leq r} S_{\mathbf{T}_j}(U) \leq \sum_{j:|j-i|\leq r} S_{\mathbf{T}_j}(U + t\delta_{\mathbf{T}_i})$$

for all $|t| \leq t_0$. Thus

$$\sum_{\mathbf{j} \in A_{i-r}} \partial_{\mathbf{j}} S_{\mathbf{T}_i - r\mathbf{e}_1}(U) + \cdots + \sum_{\mathbf{j} \in A_i} \partial_{\mathbf{j}} S_{\mathbf{T}_i}(U) + \cdots + \sum_{\mathbf{j} \in A_{i+r}} \partial_{\mathbf{j}} S_{\mathbf{T}_i + r\mathbf{e}_1}(U) = 0, \tag{3.27}$$

where

$$A_{i+k} = \{\mathbf{j} \in \mathbb{Z}^n \mid \|\mathbf{j} - \mathbf{T}_i\| \leq r - |k|, \mathbf{j}_1 = i\}$$

for $k = -r, \cdots, r$. Notice that $U \in \mathbb{R}^{\mathbb{Z} \times (\mathbb{Z}/\{1\})^{n-1}}$, thus for $k = -r, \cdots, 0, \cdots, r$, and $\mathbf{j} \in A_{i+k}$, we obtain

$$\partial_{\mathbf{j}} S_{\mathbf{T}_i + k\mathbf{e}_1}(U)$$
$$= \partial_{\mathbf{j} - \mathbf{T}_i - k\mathbf{e}_1} s(\tau^1_{-k} U|_{B^r_0})$$
$$= \partial_{\mathbf{j} - \mathbf{T}_i - k\mathbf{e}_1} s(\tau^n_{\mathbf{j}_n} \cdots \tau^2_{\mathbf{j}_2} \tau^1_{-k} U|_{B^r_0})$$
$$= \partial_{\mathbf{T}_i} S_{2\mathbf{T}_i - \mathbf{j} + k\mathbf{e}_1}(U).$$

So (3.27) implies

$$\sum_{\mathbf{j}:\|\mathbf{j}-\mathbf{T}_i\|\leq r} \partial_{\mathbf{T}_i} S_{\mathbf{j}}(U) = 0.$$

For the last assertion, using the periodic condition in $\hat{\Gamma}_1(v,w)$, we can define $J_{1,i}$, $J_{1;p,q}$ and $J_1$ along the axis parallel to axis $\mathbf{i}_1$. To be more precise, we define $J_{1,\mathbf{i}_1}^{\mathbf{i}_2,\cdots,\mathbf{i}_n}(u) = S_{\mathbf{i}}(u) - c_0$, $J_{1;p,q}^{\mathbf{i}_2,\cdots,\mathbf{i}_n}(u) = \sum_{\mathbf{i}_1=p}^{q} J_{1,\mathbf{i}_1}^{\mathbf{i}_2,\cdots,\mathbf{i}_n}(u)$ and $J_1^{\mathbf{i}_2,\cdots,\mathbf{i}_n}(u) = \liminf_{\substack{p \to -\infty \\ q \to \infty}} J_{1;p,q}^{\mathbf{i}_2,\cdots,\mathbf{i}_n}(u)$. The above proof can be modified slightly to prove the final assertion of Proposition 3.8. □

**Remark 3.9.** *If $u_k + t\delta_{\mathbf{T}_i} \in \mathcal{Y}$, (3.26) is obvious by taking $\epsilon_k \equiv 0$. However, $u_k + t\delta_{\mathbf{T}_i} \in \mathcal{Y}$ may not hold. For our choices of $\mathcal{Y}$, we can obtain (3.26) via a truncation method. Please see the proof of (A) of Theorem 3.13 (p. 15) for a typical example.*

The next result is very useful for comparison arguments. For $v \in \mathcal{M}_0$, set
$$\Gamma_1(v) = \{u \in \hat{\Gamma}_1(v-1, v+1) \mid \|u - v\|_{\mathbf{T}_i} \to 0, |i| \to \infty\}.$$

**Remark 3.10.** *For $\Gamma_1(v)$, under the same assumptions of Proposition 3.4, one can verify that (3.10), (3.11) hold and (3.8) is valid for $|i| \to \infty$.*

Set
$$c_1(v) = \inf_{u \in \Gamma_1(v)} J_1(u)$$
and
$$\mathcal{M}_1(v) = \{u \in \Gamma_1(v) \mid J_1(u) = c_1(v)\}.$$

**Theorem 3.11.** *If $s \in C^2(\mathbb{R}^{B_0^r}, \mathbb{R})$ satisfies (S1)-(S3), then $c_1(v) = 0$ and $\mathcal{M}_1(v) = \{v\}$.*

**Proof:** Noticing $v \in \Gamma_1(v)$ and $J_1(v) = 0$, we have $c_1(v) \leq 0$. To prove the reverse inequality, suppose $u \in \Gamma_1(v)$ and $J_1(u) < \infty$. Define
$$\chi_p = \begin{cases} u, & -p+1 \leq \mathbf{i}_1 \leq p-1, \\ v, & \text{otherwise.} \end{cases}$$
Then $\chi_p \in \Gamma_1(v)$. Set $\phi_p = \chi_p|_{[-p-r, p+r] \times \mathbb{Z}^{n-1}}$ and extend it as a $(2p+2r+1)$-periodic function of $\mathbf{i}_1$. Then $\phi_p \in \Gamma_0(\mathbf{l})$ with $\mathbf{l} = (2p+2r+1, 1, \cdots, 1)$, so by Proposition 3.1,

(3.28) $$0 \leq J_{1; -p-r, p+r}(\phi_p) = J_{1; -p-r, p+r}(\chi_p) = J_1(\chi_p).$$

We have
$$\begin{aligned} J_1(\chi_p) &= J_1(u) + J_{1; -p-r+1, -p+r}(\chi_p) - J_{1; -p-r+1, -p+r}(u) \\ &\quad + J_{1; p-r, p+r-1}(\chi_p) - J_{1; p-r, p+r-1}(u) - J_{1; -\infty, -p-r}(u) - J_{1; p+r, \infty}(u) \\ &=: J_1(u) - R_p(u). \end{aligned}$$

By (3.28), we obtain
$$R_p(u) \leq J_1(u).$$
Next we will prove $R_p(u) \to 0$ as $p \to \infty$ and then complete the proof of $c_1(u) = 0$. By Remark 3.10 and Proposition 3.4, $J_{1; -\infty, -p-r}(u), J_{1; p+r, \infty}(u) \to 0$ as $p \to \infty$ and similarly
$$J_{1; -p-r+1, -p+r}(\chi_p) - J_{1; -p-r+1, -p+r}(u) \to 0, \quad J_{1; p-r, p+r-1}(\chi_p) - J_{1; p-r, p+r-1}(u) \to 0$$
as $p \to \infty$, since $\tau^1_{\pm p}\chi_p, \tau^1_{\pm p} u \to v$ in $\mathbb{R}^{\mathbb{Z}^n}$ via (3.8).

**Remark 3.12.** *In the definition of $\Gamma_1(v)$, if $v \pm 1$ is replaced by $v \pm j$ for any $j \in \mathbb{N}$, the above argument still holds.*

What is left is to show that $\mathcal{M}_1(v) = \{v\}$. Let $u \in \mathcal{M}_1(v)$. Then $v - 1 \leq u \leq v + 1$, so for any $i \in \mathbb{Z}$, the function $\delta_{\mathbf{T}_i}$ and $|t| \leq 1$, $v - 2 \leq u + t\delta_{\mathbf{T}_i} \leq v + 2$. Hence by Remark 3.12, and $u_k = u$, the assumption of Proposition 3.8 (with $\epsilon_k = 0$) is satisfied. Consequently, $u$ satisfies (2.1) for all $\mathbf{i} \in \mathbb{Z}^n$. $u \in \mathcal{M}_1(v)$ implies $\tau^1_{-1}u \in \mathcal{M}_1(v)$. If $\tau^1_{-1}u = u$, i.e., $u(\mathbf{i} + \mathbf{e}_1) = u(\mathbf{i})$, together with the fact $\|u - v\|_{\mathbf{T}_i} \to 0$ as $|i| \to \infty$, we have $u = v$, completing the proof. Thus assume $u \neq \tau^1_{-1}u$. We claim that

(3.29) $\qquad (i)\ u < \tau^1_{-1}u \quad \text{or} \quad (ii)\ u > \tau^1_{-1}u.$

Indeed, set $\phi = \max(u, \tau^1_{-1}u)$ and $\psi = \min(u, \tau^1_{-1}u)$. Noticing that for any $i \in \mathbb{Z}$, by Lemma 2.8 we have
$$S_{\mathbf{T}_i}(\phi) + S_{\mathbf{T}_i}(\psi) \leq S_{\mathbf{T}_i}(u) + S_{\mathbf{T}_i}(\tau^1_{-1}u),$$

and then

(3.30) $$J_{1,i}(\phi) + J_{1,i}(\psi) \leq J_{1,i}(u) + J_{1,i}(\tau^1_{-1}u).$$

Therefore summing over $i$ leads to

(3.31) $$J_1(\phi) + J_1(\psi) \leq J_1(u) + J_1(\tau^1_{-1}u) = 0.$$

Since $\phi, \psi \in \Gamma_1(v)$, $J_1(\phi), J_1(\psi) \geq c_1(v) = 0$. Hence by (3.31), $\phi, \psi \in \mathcal{M}_1(v)$ and thus they satisfy (2.1) (by the above argument that $u \in \mathcal{M}_1(v)$ implies $u$ satisfies (2.1)). Consequently by Corollary 2.7 we prove our claim (3.29).

Assume that (3.29) (i) holds. Then for all $j \in \mathbb{N}$,

$$\tau^1_j u < u < \tau^1_{-j} u.$$

Letting $j \to \infty$ gives

$$v \leq u \leq v.$$

Similarly one can prove the case of (3.29) (ii) and thus we complete the proof of Therorem 3.11. □

We are now in a position to state our first main result of this paper. To this end, assuming $v_0 < w_0$ are adjacent members in $M_0$, denote $\Gamma_1 = \Gamma_1(v_0, w_0)$ and define

(3.32) $$c_1 = c_1(v_0, w_0) = \inf_{u \in \Gamma_1(v_0, w_0)} J_1(u),$$

and

$$\mathcal{M}_1 = \mathcal{M}_1(v_0, w_0) = \{u \in \Gamma_1(v_0, w_0) \,|\, J_1(u) = c_1\}.$$

**Theorem 3.13.** *If $s \in C^2(\mathbb{R}^{B_0^r}, \mathbb{R})$ satisfies (S1)-(S3) and $(*_0)$ holds, then there is a solution $U_1 \in \mathcal{M}_1$ of (2.1). Moreover, $\mathcal{M}_1$ is an ordered set and the elements of $\mathcal{M}_1$ are solutions of (2.1), and any $U \in \mathcal{M}_1$ is strictly 1-monotone in $\mathbf{i}_1$.*

**Proof:** Taking $(u_k) \subset \Gamma_1$ as a minimizing sequence for (3.32), dropping finite terms if necessary, we have that $J_1(u_k) \leq M$ holds for some $M > 0$ and for all $k \in \mathbb{N}$. Note that $\Gamma_1$ is not a complete space in pointwise convergence since $v_0, w_0$ are limit points of some sequences in $\Gamma_1$. To obtain an element having the asymptotic properties of $\Gamma_1$, noticing $u \in \Gamma_1$ implies $\tau^1_{-j}u \in \Gamma_1$ for all $j \in \mathbb{Z}$, we may assume that

(3.33) $$u_k(\mathbf{T}_i) \leq \frac{1}{2}(v_0 + w_0)(\mathbf{T}_0) \leq u_k(\mathbf{T}_0)$$

for all $-i, k \in \mathbb{N}$.

By Proposition 3.7 there is a $U_1 \in \hat{\Gamma}_1(v_0, w_0)$ such that $u_k \to U_1$ pointwise along a subsequence and

(3.34) $$-K_1 \leq J_1(U_1) \leq c_1 + 1 + 2K_1.$$

Without loss of generality we may take this subsequence as $(u_k)$. By (3.33), for $0 > i \in \mathbb{Z}$,

(3.35) $$U_1(\mathbf{T}_i) \leq \frac{1}{2}(v_0 + w_0)(\mathbf{T}_0) \leq U_1(\mathbf{T}_0),$$

so $v_0 \not\equiv U_1 \not\equiv w_0$.

To complete the proof of Theorem 3.13, we will show that:

(A) $U_1$ is a solution of (2.1), as is any $U \in \mathcal{M}_1$;
(B) $U_1$ and any $U \in \mathcal{M}_1$ are strictly 1-monotone in $\mathbf{i}_1$;
(C) $J_1(U_1) = c_1$, so $\mathcal{M}_1 \neq \emptyset$;
(D) $\mathcal{M}_1$ is an ordered set.

*Proof of* (A). To prove the first assertion, we only need to show that the assumption of Proposition 3.8 for $(u_k)$ holds. Since $v_0 \leq u_k + t\delta_{\mathbf{T}_i} \leq w_0$ may be not hold, we use a truncation trick to recover this condition. Since $v_0 \leq u_k \leq w_0$, for $|t_0| \leq 1$,

$$w_0 - 2 \leq v_0 - 1 \leq u_k + t\delta_{\mathbf{T}_i} \leq w_0 + 1 \leq w_0 + 2. \tag{3.36}$$

Define $f_k = \max(u_k + t\delta_{\mathbf{T}_i}, w_0)$ and $g_k = \min(u_k + t\delta_{\mathbf{T}_i}, w_0)$. Of course $f_k \in \Gamma_1(w_0)$. Thus by Theorem 3.11,

$$J_1(f_k) \geq 0, \tag{3.37}$$

and thus

$$J_1(g_k) \leq J_1(f_k) + J_1(g_k). \tag{3.38}$$

Noting $g_k \in \hat{\Gamma}_1(v_0 - 1, w_0)$ and proceeding as in (3.30)-(3.31) gives

$$J_1(f_k) + J_1(g_k) \leq J_1(u_k + t\delta_{\mathbf{T}_i}) + J_1(w_0) = J_1(u_k + t\delta_{\mathbf{T}_i}). \tag{3.39}$$

Define $\chi_k = \max(g_k, v_0)$ and $\psi_k = \min(g_k, v_0)$. Then $\chi_k \in \Gamma_1$ and $\psi_k \in \Gamma(v_0)$, so as in (3.37)-(3.39),

$$J_1(\chi_k) \leq J_1(\chi_k) + J_1(\psi_k) \leq J_1(g_k) + J_1(v_0) = J_1(g_k). \tag{3.40}$$

Combining (3.38)-(3.40) gives

$$c_1 \leq J_1(u_k) =: c_1 + \epsilon_k \leq J_1(\chi_k) + \epsilon_k \leq J_1(u_k + t\delta_{\mathbf{T}_i}) + \epsilon_k,$$

where $\epsilon_k \to 0$ as $k \to \infty$. Thus the assumption of Proposition 3.8 holds and $U_1$ is a solution of (2.1). Next for $U \in \mathcal{M}_1$, we see the sequence $\{\phi_k \,|\, \phi_k \equiv U\}$ is a minimizing sequence for (3.32). Hence we can now proceed analogously to the above proof to show that $U$ is a solution of (2.1).

*Proof of* (B). Assuming $U_1$ is 1-monotone in $\mathbf{i}_1$, i.e.,

$$U_1 \leq \tau^1_{-1} U_1, \tag{3.41}$$

and noticing $U_1 \in \hat{\Gamma}_1(v_0, w_0) \setminus \{v_0, w_0\}$, by (3.34), we obtain that Corollary 3.6 implies $U_1 \in \Gamma_1(v_0, w_0)$. Similarly, any $U \in \mathcal{M}_1$ belongs to $\Gamma_1(v_0, w_0)$.

Now we prove (3.41). To this end, define $\Phi_k = \max(u_k, \tau^1_{-1} u_k)$ and $\Psi_k = \min(u_k, \tau^1_{-1} u_k)$. Then $\Phi_k, \Psi_k \in \Gamma_1$ and as in (3.30)-(3.31),

$$J_1(\Phi_k) + J_1(\Psi_k) \leq J_1(u_k) + J_1(\tau^1_{-1} u_k) = 2J_1(u_k). \tag{3.42}$$

Since $J_1(u_k) \to c_1$ as $k \to \infty$, we obtain that $\Phi_k$ and $\Psi_k$ are minimizing sequences for (3.32). By Propositions 3.7 and 3.8, and noticing $\max(\cdot, \cdot)$ and $\min(\cdot, \cdot)$ are continuous on $\mathbb{R}^{\mathbb{Z}^n}$ we have

$$\Phi_k \to \Phi = \max(U_1, \tau^1_{-1} U_1) \quad \text{and} \quad \Psi_k \to \Psi = \min(U_1, \tau^1_{-1} U_1)$$

as $k \to \infty$. Thus $\Phi, \Psi$ are solutions of (2.1). By Corollary 2.7, we have (a) $U_1 \equiv \tau^1_1 U_1$ or (b) $U_1 > \tau^1_{-1} U_1$ or (c) $U_1 < \tau^1_{-1} U_1$. If (a) occurs, $U_1$ is a constant along $\mathbf{i}_1$, so $U_1 \in \Gamma_0 \cap \hat{\Gamma}_1$. Moreover, note that $v_0 \not\equiv U_1 \not\equiv w_0$. Therefore $J_0(U_1) > c_0$, so $J_1(U_1) = \infty$, contrary to (3.34). Notice that (b) is contrary to the asymptotic properties in $\Gamma_1$ or to (3.35) (for $i = -1$). So (c) holds and (3.41) is proved.

The above proof can be applied to prove that any $U \in \mathcal{M}_1$ satisfies (3.41). Indeed, replacing $U_1, u_k$ in the above paragraph by $U$ and noting $J_1(U) = c_1$, which is used in (3.42), we have $U \in \mathcal{M}_1$ also satisfies (3.41). The strict inequalities of (3.41) for $U_1$ or $U \in \mathcal{M}_1$ follow from Lemma 2.6.

*Proof of* (C). Firstly (3.8)-(3.9) imply

(3.43) $$\begin{cases} \|U_1 - v_0\|_{\mathbf{T}_i} \to 0, & i \to -\infty, \\ \|U_1 - w_0\|_{\mathbf{T}_i} \to 0, & i \to \infty. \end{cases}$$

Define $\hat{\mathbf{T}}_i = \cup_{j=i-r}^{i+r} \mathbf{T}_j$. For any $\epsilon > 0$, by (3.43), there is a $p_0 = p_0(\epsilon)$ such that if $p \geq p_0$,

(3.44) $$\|U_1 - v_0\|_{\hat{\mathbf{T}}_{-p}} \leq \epsilon/2, \quad \|U_1 - w_0\|_{\hat{\mathbf{T}}_p} \leq \epsilon/2.$$

Since $u_k \to U_1$ as $k \to \infty$, then for any $p \geq p_0$, there is a $k_0 = k_0(p)$ such that for $k \geq k_0$,

(3.45) $$\|u_k - U_1\|_{\hat{\mathbf{T}}_{-p}} \leq \epsilon/2, \quad \|u_k - U_1\|_{\hat{\mathbf{T}}_p} \leq \epsilon/2.$$

Thus for such $k$ and $p$,

(3.46) $$\|u_k - v_0\|_{\hat{\mathbf{T}}_{-p}} \leq \epsilon, \quad \|u_k - w_0\|_{\hat{\mathbf{T}}_p} \leq \epsilon.$$

For fixed $k \geq k_0(p)$, noting $u_k \in \Gamma_1$, we obtain a $q_0 = q_0(k)$ such that for $q \geq q_0$,

(3.47) $$\|u_k - v_0\|_{\hat{\mathbf{T}}_{-q}} \leq \epsilon, \quad \|u_k - w_0\|_{\hat{\mathbf{T}}_q} \leq \epsilon.$$

Define

(3.48) $$f_k = \begin{cases} w_0, & p - r \leq \mathbf{i}_1 \leq p + r \quad \text{or} \quad q - r \leq \mathbf{i}_1 \leq q + r, \\ u_k, & p + r + 1 \leq \mathbf{i}_1 \leq q - r - 1, \end{cases}$$

and

(3.49) $$g_k = \begin{cases} v_0, & -q - r \leq \mathbf{i}_1 \leq -q + r \quad \text{or} \quad -p - r \leq \mathbf{i}_1 \leq -p + r, \\ u_k, & -q + r + 1 \leq \mathbf{i}_1 \leq -p - r - 1. \end{cases}$$

Extend $f_k$ (resp. $g_k$) to a $(q + 2r + 1 - p)$-periodic function of $\mathbf{i}_1$ and still denote it by $f_k$ (resp. $g_k$). Then by (3.46)-(3.47), there is a $\kappa(\epsilon)$ such that

(3.50) $$\begin{aligned} |J_{1;p,q}(u_k) - J_{1;p,q}(f_k)| &\leq \kappa(\epsilon), \\ |J_{1;-q,-p}(u_k) - J_{1;-q,-p}(g_k)| &\leq \kappa(\epsilon) \end{aligned}$$

and $\kappa(\epsilon) \to 0$ as $\epsilon \to 0$. By Proposition 3.1,

(3.51) $$\begin{aligned} J_{1;p,q}(f_k) &= J_{1;p-r,q+r}(f_k) \geq 0, \\ J_{1;-q,-p}(g_k) &= J_{1;-q-r,-p+r}(g_k) \geq 0. \end{aligned}$$

Since

(3.52) $$\begin{aligned} J_{1;1,\infty}(u_k) &= J_{1;1,p-1}(u_k) + J_{1;p,q}(u_k) + J_{1;q+1,\infty}(u_k), \\ J_{1;-\infty,0}(u_k) &= J_{1;-\infty,-q-1}(u_k) + J_{1;-q,-p}(u_k) + J_{1;-p+1,0}(u_k), \end{aligned}$$

by (3.50)-(3.52),

(3.53) $$\begin{aligned} J_{1;1,\infty}(u_k) &\geq J_{1;1,p-1}(u_k) - \kappa(\epsilon) + J_{1;q+1,\infty}(u_k), \\ J_{1;-\infty,0}(u_k) &\geq J_{1;-\infty,-q-1}(u_k) - \kappa(\epsilon) + J_{1;-p+1,0}(u_k). \end{aligned}$$

Letting $q \to \infty$ in (3.53) gives

(3.54) $$J_1(u_k) \geq J_{1;-p+1,p-1}(u_k) - 2\kappa(\epsilon).$$

Thus letting $k \to \infty$ shows that

(3.55) $$c_1 \geq J_{1;-p+1,p-1}(U_1) - 2\kappa(\epsilon).$$

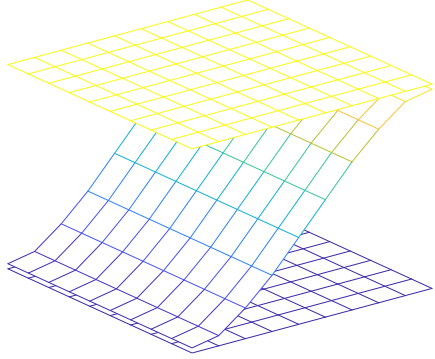

FIGURE 1. Elements in $\mathcal{M}_0 \cup \mathcal{M}_1(v_0, w_0)$.

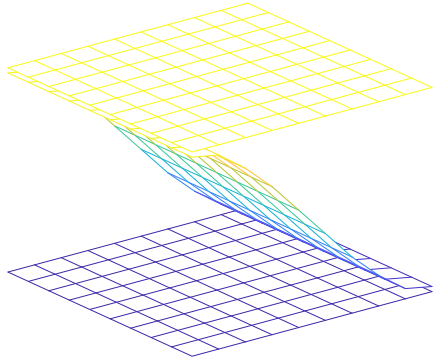

FIGURE 2. Elements in $\mathcal{M}_0 \cup \mathcal{M}_1(w_0, v_0)$.

Lastly, letting $p \to \infty$ and then $\epsilon \to 0$ yields

$$c_1 \geq J_1(U_1). \tag{3.56}$$

The reverse inequality follows from $U_1 \in \Gamma_1$ and thus the proof of the first assertion of Theorem 3.13 is complete.

*Proof of* (D). Let $V, W \in \mathcal{M}_1$. Define $\Phi = \max(V, W)$ and $\Psi = \min(V, W)$ and then we have

$$J_1(\Phi) + J_1(\Psi) \leq J_1(V) + J_1(W) = 2c_1. \tag{3.57}$$

Proceeding as in the proof of (B), $\Phi, \Psi \in \mathcal{M}_1$. By (A) and Corollary 2.7 we complete the proof of (D) and thus Theorem 3.13. □

Figure 1 shows a typical configuration in $\mathcal{M}_1$.

**Remark 3.14.** *Similar to $\mathcal{M}_1(v_0, w_0)$, one can define $\mathcal{M}_1(w_0, v_0)$ as follows:*

$$\mathcal{M}_1(w_0, v_0) = \{u \in \Gamma_1(w_0, v_0) \,|\, J_1(u) = c_1(w_0, v_0)\},$$

*where*
$$\Gamma_1(w_0, v_0) = \{u \in \hat{\Gamma}_1(v_0, w_0) \mid \|u - w_0\|_{\mathbf{T}_i} \to 0, i \to -\infty; \|u - v_0\|_{\mathbf{T}_i} \to 0, i \to \infty\}.$$
*and*
$$c_1(w_0, v_0) = \inf_{u \in \Gamma_1(w_0, v_0)} J_1(u).$$

$\mathcal{M}_1(w_0, v_0)$ *is a nonempty and ordered set and the elements in* $\mathcal{M}_1(w_0, v_0)$ *have analogous properties as in Theorem 3.13. The proof is the same to that of Theorem 3.13. See Figure 2 for an example of configurations of* $\mathcal{M}_1(w_0, v_0)$.

**Remark 3.15.** *Assume* $v, w \in \mathcal{M}_0$. *By Theorem 3.13, if* $v, w$ *are adjacent then* $\mathcal{M}_1(v, w) \neq \emptyset$. *Conversely, if* $\mathcal{M}_1(v, w) \neq \emptyset$, *i.e., there are minimal and Birkhoff configuration heteroclinic to* $v$ *and* $w$ *in* $i_1$, *then* $v, w$ *are adjacent members in* $\mathcal{M}_0$. *For the proof of this fact, we refer the reader to [22, Theorem 3.34].*

The critical value $c_1$ can be characterized in another way. To this end, set
$$\mathcal{S}_1 = \{u \in \hat{\Gamma}_1(v_0, w_0) \mid u \leq \tau^1_{-1} u \text{ and } v_0 \not\equiv u \not\equiv w_0\},$$
then we have $c_1 = \inf_{u \in \mathcal{S}_1} J_1(u)$. The proof is easy (cf. [22, Corollary 3.32]) and we omit it.

If we perturb $s$ slightly, the gap condition $(*_0)$ still holds. To state this fact more precisely, we need some notations. Suppose $\bar{s} \in C^2(\mathbb{R}^{B_0^r}, \mathbb{R})$ satisfies (S1)-(S3). For $u \in \Gamma_0$, set
$$S_{\mathbf{j}}^{\bar{s}}(u) = \bar{s}(\tau^n_{-\mathbf{j}_n} \cdots \tau^1_{-\mathbf{j}_1} u|_{B_0^r}),$$
$$J_0^{\bar{s}}(u) = S_0^{\bar{s}}(u);$$
$$c_0(\bar{s}) = \inf_{u \in \Gamma_0} J_0^{\bar{s}}(u);$$
and
$$\mathcal{M}_0(\bar{s}) = \{u \in \Gamma_0 \mid J_0^{\bar{s}}(u) = c_0(\bar{s})\}.$$
When $(*_0)$ holds for $\bar{s}$, we denote the associated gap pair by $v_0(\bar{s}), w_0(\bar{s})$.

**Proposition 3.16.** *Assume that* $s \in C^2(\mathbb{R}^{B_0^r}, \mathbb{R})$ *satisfies (S1)-(S3) and* $(*_0)$ *holds for* $s$. *There is an* $\epsilon$ *such that for* $\bar{s} \in C^2(\mathbb{R}^{B_0^r}, \mathbb{R})$ *satisfies (S1)-(S3), if*

(3.58) $$\|s - \bar{s}\|_{L^\infty(\mathbb{R}^{\#B_0^r})} \leq \epsilon,$$

*then* $(*_0)$ *holds for* $\bar{s}$. *Moreover, suppose* $v_0$, $w_0$ *are a gap pair for* $s$ *and*
$$\alpha_0 = v_0(\mathbf{0}); \quad \beta_0 = w_0(\mathbf{0}).$$
*Then for any* $\delta \in (0, \frac{\beta_0 - \alpha_0}{2})$, *there is an* $\epsilon_1 = \epsilon_1(s, \delta)$ *such that (3.58) holds with* $\epsilon = \epsilon_1$ *and*

(3.59) $$v(\mathbf{0}) \notin (\alpha_0 + \delta, \beta_0 - \delta)$$

*for all* $v \in \mathcal{M}_0(\bar{s})$.

**Proof:** Clearly the second assertion implies the first one, so we only prove the second one. We prove it by contradiction arguments. If it is not true, for some $\delta$ we have a sequence $(s_k)$ satisfying (S1)-(S3),

(3.60) $$\|s - s_k\|_{L^\infty(\mathbb{R}^{\#B_0^r})} \leq \frac{1}{k},$$

and an associated $u_k \in \mathcal{M}_0(s_k)$ with

(3.61) $$u_k(\mathbf{0}) \in (\alpha_0 + \delta, \beta_0 - \delta).$$

By the periodicity of $u_k$, we can assume that there is a $u \in \Gamma_0(s)$ such that $u_k \to u$ pointwise (taking a subsequence if necessary). By (3.60), $J_0(u) = c_0(s)$, so $u \in \mathcal{M}_0(s)$. But by (3.61),
$$u(\mathbf{0}) \in (\alpha_0, \beta_0),$$
contrary to $(*_0)$ for $s$. $\square$

**Remark 3.17.** *For a gap pair $v_0(s), w_0(s)$ for $s$, the proof of Proposition 3.16 shows that the unique gap pair $v_0(\bar{s}), w_0(\bar{s})$ for $(*_0)$ for $\bar{s}$ approaches $v_0(s), w_0(s)$ as $\bar{s} \to s$ in $C(\mathbb{R}^{B_\mathbf{0}^r}, \mathbb{R})$.*

The condition $(*_0)$ is generic as we see from the next proposition. This is because of our choice of rotation vector $\alpha = \mathbf{0}$, which is rational.

**Proposition 3.18.** *Assume that $s \in C^2(\mathbb{R}^{B_\mathbf{0}^r}, \mathbb{R})$ satisfies (S1)-(S3). Then there is an $\bar{s} \in C^2(\mathbb{R}^{B_\mathbf{0}^r}, \mathbb{R})$ satisfying (S1)-(S3) such that if $\epsilon > 0$, $(*_0)$ holds for (2.1) with $s$ replaced by $s + \epsilon\bar{s}$.*

**Proof:** Let $v \in \mathcal{M}_0(s)$ and set
$$\bar{s}(u) = \sin^2 \pi(u(\mathbf{0}) - v(\mathbf{0})) + \frac{1}{2} \sum_{\mathbf{k}: \|\mathbf{k}\| \leq r} (u(\mathbf{k}) - u(\mathbf{0}))^2.$$
Then $\bar{s}$ satisfies all the conditions. Indeed, $\bar{s}$ satisfies (S1)-(S3) and for any $\epsilon > 0$,
$$\mathcal{M}_0(s + \epsilon\bar{s}) = \{v(\mathbf{0}) + j \mid j \in \mathbb{Z}\}.$$
$\square$

**Remark 3.19.** *Proposition 3.18 can be stated in another form.*

*If $s \in C^2(\mathbb{R}^{B_\mathbf{0}^r}, \mathbb{R})$ satisfies (S1)-(S3). Then for any $\epsilon > 0$, there is an $\bar{s} \in C^2(\mathbb{R}^{B_\mathbf{0}^r}, \mathbb{R})$ satisfying (S1)-(S3) such that $(*_0)$ holds for (2.1) with $s$ replaced by $\bar{s}$, and*
$$\|s - \bar{s}\|_{L^\infty(\mathbb{R}^{\#B_\mathbf{0}^r})} + \sum_{\mathbf{j} \in B_\mathbf{0}^r} \|\partial_\mathbf{j} s - \partial_\mathbf{j} s_k\|_{L^\infty(\mathbb{R}^{\#B_\mathbf{0}^r})} \leq \epsilon.$$
*Indeed, $\bar{s}(u) := s(u) + \frac{1}{4\pi}\epsilon \sin^2 \pi(u(\mathbf{0}) - v(\mathbf{0}))$ satisfies all the conditions.*

Proposition 3.1 shows that $\mathcal{M}_0 = \mathcal{M}_0(\mathbf{1})$. Maybe someone wants to search other periodic solutions with period other than 1. But the following proposition tells us there is no such solutions. Extend the definition of $J_1(u)$ to $J_1(\mathbf{l}, u)$ etc. We have

**Proposition 3.20.** $\mathcal{M}_1(\mathbf{l}) = \mathcal{M}_1$ *and* $c_1(\mathbf{l}) := \inf_{u \in \Gamma_1(\mathbf{l})} J_1(\mathbf{l}, u) = (\prod_{i=2}^n \mathbf{l}_i) c_1$.

**Proof:** In fact, if we replace
$$u(\cdot + \mathbf{e}_i) = u(\cdot),\ 2 \leq i \leq n$$
by
$$(3.62) \qquad\qquad u(\cdot + \mathbf{l}_i \mathbf{e}_i) = u(\cdot),\ 2 \leq i \leq n$$
for some $\mathbf{l} = (\mathbf{l}_2, \cdots, \mathbf{l}_n) \in \mathbb{N}^{n-1}$, the results of this section are also true. For example, the conclusion of Theorem 3.13 holds. Thus $\mathcal{M}_1(\mathbf{l})$ is ordered. Using this fact and $u(\cdot + \mathbf{e}_i) \in \mathcal{M}_1(\mathbf{l})$, $2 \leq i \leq n$, we can prove Proposition 3.20 exactly as in Proposition 3.1. $\square$

The last theorem of this section explores the relation of solutions of (2.1) that are minimal and Birkhoff and the solutions of (2.1) in $\mathcal{M}_0$, $\mathcal{M}_1(v_0, w_0)$ and $\mathcal{M}_1(w_0, v_0)$. Note that in the assumption of (2) of Theorem 3.21, $u$ is heteroclinic in $\mathbf{i}_1$. This assumption ensures the results similar to Aubry-Mather theory and if $u$ is heteroclinic in more directions, the case becomes complex. Theorem 3.21 also clarify part of the relations of our heteroclinic solutions with the solutions obtained by Miao, et al. [14].

**Theorem 3.21.** *Assume that $s \in C^2(\mathbb{R}^{B_\mathbf{0}^r}, \mathbb{R})$ satisfies (S1)-(S3).*

(1) *If $u \in \mathcal{M}_0$ or if $(*_0)$ holds and $u \in \mathcal{M}_1(v_0, w_0) \cup \mathcal{M}_1(w_0, v_0)$, then $u$ is minimal and Birkhoff.*
(2) *If $u \in \mathbb{R}^{\mathbb{Z} \times (\mathbb{Z}/\{1\})^{n-1}}$ is a minimal and Birkhoff solution of (2.1) with rotation vector $\mathbf{0}$, then $u \in \mathcal{M}_0$ or $(*_0)$ holds and $u \in \mathcal{M}_1(v_0, w_0) \cup \mathcal{M}_1(w_0, v_0)$ for some adjacent pair $v_0, w_0 \in \mathcal{M}_0$.*

**Proof:** (1) If $u \in \mathcal{M}_0$ then it is Birkhoff following from Theorem 2.12. By [16, Theorem 4.8], $u$ is minimal. Now suppose $u \in \mathcal{M}_1(v_0, w_0) \cup \mathcal{M}_1(w_0, v_0)$. By Theorem 3.13, $u$ is Birkhoff. To prove $u$ is minimal, just proceed as in the proof of [16, Theorem 4.8] and use Proposition 3.20 to obtain a contradiction. Thus we complete the proof of (1) of Theorem 3.21.

To prove (2), we need the following technical lemma.

**Lemma 3.22.** *If $u \in \mathbb{R}^{\mathbb{Z} \times (\mathbb{Z}/\{1\})^{n-1}}$ is minimal, then for any $\phi \in \mathbb{R}^{\mathbb{Z} \times (\mathbb{Z}/\{1\})^{n-1}}$ with compact support in $\mathbf{i}_1$,*

$$\sum_{\mathbf{j} \in \mathbb{Z} \times \{0\}^{n-1}} S_\mathbf{j}(u + \phi) - S_\mathbf{j}(u) \geq 0. \tag{3.63}$$

*Proof of Lemma 3.22.* Define

$$\theta_l(i) = \begin{cases} 1, & |i| \leq l, \\ 0, & |i| > l \end{cases}$$

and $(\theta \psi)(\mathbf{j}) = \theta_l(|\mathbf{j}_2|) \cdots \theta_l(|\mathbf{j}_n|) \psi(|\mathbf{j}|)$. Since $u$ is minimal,

$$0 \leq \sum_{\mathbf{j} \in \mathbb{Z}^n} S_\mathbf{j}(u + \theta \phi) - S_\mathbf{j}(u). \tag{3.64}$$

Suppose the support of $\phi$ lies in $[p, q] \times \mathbb{Z}^{n-1}$; with $p, q \in \mathbb{Z}$. Then by (3.64), we have

$$\begin{aligned}
0 &\leq \sum_{\mathbf{j} \in [p-r, q+r] \times [-l-r, l+r]^{n-1}} [S_\mathbf{j}(u + \theta \phi) - S_\mathbf{j}(u)] \\
&= \sum_{\mathbf{j} \in [p-r, q+r] \times [-l+r, l-r]^{n-1}} [S_\mathbf{j}(u + \phi) - S_\mathbf{j}(u)] + \mathcal{R}_l(u, \phi) \\
&= (2l - 2r + 1)^{n-1} \sum_{\mathbf{j} \in [p-r, q+r] \times \{0\}^{n-1}} [S_\mathbf{j}(u + \phi) - S_\mathbf{j}(u)] + \mathcal{R}_l(u, \phi),
\end{aligned} \tag{3.65}$$

where

$$\mathcal{R}_l(u, \phi) = \sum_{A_l} [S_\mathbf{j}(u + \theta \phi) - S_\mathbf{j}(u)]$$

and $A_l$ is the region

$$[p - r, q + r] \times ([-l - r, l + r]^{n-1} \setminus [-l + r, l - r]^{n-1}).$$

By Lemma 2.11,

$$\sum_{A_l} [S_{\mathbf{j}}(u + \theta\phi) - S_{\mathbf{j}}(u)]$$

(3.66)
$$\leq L \sum_{\bar{A}_l} |\theta_l \phi|$$

$$\leq L \sum_{\bar{A}_l} |\phi|$$

$$\leq L(\#B_{\mathbf{0}}^r)(q + 1 - p)[(2l + 2r + 1)^{n-1} - (2l - 2r + 1)^{n-1}]M,$$

where $M = \sup_{\mathbf{j} \in [p,q] \times \{0\}^{n-1}} |\phi(\mathbf{j})|$. Therefore by (3.65)-(3.66),

(3.67)
$$0 \leq (2l - 2r + 1)^{n-1} \sum_{\mathbb{Z} \times \{0\}^{n-1}} [S_{\mathbf{j}}(u + \phi) - S_{\mathbf{j}}(u)]$$
$$+ L(\#B_{\mathbf{0}}^r)(q + 1 - p)[(2l + 2r + 1)^{n-1} - (2l - 2r + 1)^{n-1}]M.$$

Dividing (3.67) by $(2l - 2r + 1)^{n-1}$ and letting $l \to \infty$ yields (3.63). This proves Lemma 3.22.

*Proof of (2) of Theorem 3.21.* Let $u \in \mathbb{R}^{\mathbb{Z} \times (\mathbb{Z}/\{1\})^{n-1}}$ be a minimal and Birkhoff solution of (2.1) with rotation vector $\alpha = \mathbf{0}$. Since $\alpha = \mathbf{0}$, by Lemma 2.5, $u$ is bounded. Thus we have a smallest $w$ and largest $v$ in $\mathcal{M}_0$ such that $v \leq u \leq w$. If $v = w$, $u \in \mathcal{M}_0$. Therefore suppose that $v < w$ and $u \neq v, w$. By Lemma 2.6, $v(\mathbf{i}) < u(\mathbf{i}) < w(\mathbf{i})$ for all $\mathbf{i}$. Since $u$ is Birkhoff, $\tau_{-1}^1 u = u$, $\tau_{-1}^1 u > u$, or $\tau_{-1}^1 u < u$. The cases of $\tau_{-1}^1 u > u$ and $\tau_{-1}^1 u < u$ are treated similarly, so we only prove the case of

(3.68)
$$\tau_{-1}^1 u < u.$$

Define $u_k = \tau_k^1 u$ for $k \in \mathbb{Z}$. Since $v \leq u \leq w$, we have $v \leq u_k \leq w$ and by (3.68),

(3.69)
$$u_{k+1} > u_k.$$

Thus $u_k$ converges to $\bar{u} \leq w$ (resp. $\underline{u} \geq v$) as $k \to \infty$ (resp. $k \to -\infty$). (3.69) implies $\tau_{-1}^1 \bar{u} = \bar{u}$ and $\tau_{-1}^1 \underline{u} = \underline{u}$. Thus $\underline{u}, \bar{u} \in \Gamma_0$.

We claim that $\underline{u}, \bar{u} \in \mathcal{M}_0$. Indeed, if

(3.70)
$$J_0(\bar{u}) > c_0,$$

then there is a $k_0 \in \mathbb{N}$ such that for $k \geq k_0$,

(3.71)
$$J_0(u_k) - c_0 \geq \frac{1}{2}(J_0(\bar{u}) - c_0) =: \gamma > 0.$$

Therefore for $q \geq p + 2r + 1 \geq p - r \geq k_0$,

(3.72)
$$J_{1;p,q}(u) \geq (q + 1 - p)\gamma.$$

Set

$$f_{p,q} = \begin{cases} w, & p + 1 \leq \mathbf{i}_1 \leq q - 1, \\ u, & \text{otherwise.} \end{cases}$$

By Lemma 3.22,

(3.73)
$$0 \geq J_{1;p-r,q+r}(u) - J_{1;p-r,q+r}(f_{p,q}).$$

Then by (3.72)-(3.73),

$$0 \geq J_{1;p-r,q+r}(u) - J_{1;p-r,q+r}(f_{p,q})$$

(3.74)
$$= [\sum_{i=p-r}^{p+r} + \sum_{i=p+r+1}^{q-r-1} + \sum_{i=q-r}^{q+r}][J_{1,i}(u) - J_{1,i}(f_{p,q})]$$

$$\geq (q-p-2r-1)\gamma + [\sum_{i=p-r}^{p+r} + \sum_{i=q-r}^{q+r}][J_{1,i}(u) - J_{1,i}(f_{p,q})]$$

By Lemma 2.11, it is easy to see that the last two terms on the right in (3.74) are bounded. Hence (3.74) cannot hold for $q - p$ sufficiently large. Thus $J_0(\overline{u}) = c_0$. Similarly we have $J_0(\underline{u}) = c_0$. Thus $\underline{u} = v$, $\overline{u} = w$, and $u \in \Gamma_1(w, v)$. If $\tau_{-1}^1 u = u$, the above argument shows $u \in \mathcal{M}_0$, which contradicts the choices of $v, w$.

Next we claim

(3.75)
$$J_1(u) = c_1(w, v).$$

If the claim holds, by Remark 3.15, $v, w \in \mathcal{M}_0$ are adjacent. So $(*_0)$ holds and $u \in \mathcal{M}_1(w, v)$. If (3.75) is false, since $u \in \Gamma_1(w, v)$,

(3.76)
$$J_1(u) > c_1(w, v).$$

To exclude the case of $J_1(u) = \infty$, we choose $U \in \Gamma_1(w, v)$ such that for some $\sigma > 0$,

(3.77)
$$c_1 \leq J_1(U) < J_1(U) + \sigma < J_1(u).$$

For any $\kappa > 0$, there is a $q = q(\kappa) \in \mathbb{N}$ such that for $\phi \in \{u, U\}$,

(3.78)
$$\begin{cases} \|\phi - w\|_{\mathbf{T}_i} \leq \kappa & i \leq -q, \\ \|\phi - v\|_{\mathbf{T}_i} \leq \kappa & i \geq q. \end{cases}$$

For $p \in \mathbb{N}$ and $p > q + 2r + 2$, set

$$\psi = \begin{cases} U, & -p+1 \leq \mathbf{i}_1 \leq p-1, \\ u, & \text{otherwise.} \end{cases}$$

Thus for $\kappa = \kappa(\sigma)$ sufficiently small and $\phi \in \{u, U, \psi\}$,

(3.79)
$$|S_\mathbf{i}(\phi)| \leq \frac{\sigma}{12(r+1)}$$

for $|\mathbf{i}_1| > q(\kappa) + r$. For $p$ sufficiently large,

(3.80)
$$J_{1;-p,p}(U) \leq J_1(U) + \sigma/6.$$

Thus we have

$$\sum_{\mathbf{j}\in[-p-r,p+r]\times\{0\}^{n-1}} (S_{\mathbf{j}}(u) - S_{\mathbf{j}}(U))$$

$$= \sum_{\mathbf{j}\in[-p-r,p+r]\times\{0\}^{n-1}} (S_{\mathbf{j}}(u) - S_{\mathbf{j}}(\psi)) + \sum_{\mathbf{j}\in[-p-r,p+r]\times\{0\}^{n-1}} (S_{\mathbf{j}}(\psi) - S_{\mathbf{j}}(U))$$

$$= \sum_{\mathbf{j}\in[-p-r,p+r]\times\{0\}^{n-1}} (S_{\mathbf{j}}(u) - S_{\mathbf{j}}(\psi))$$

$$+ \left[\sum_{\mathbf{j}\in[-p-r,-p+r]\times\{0\}^{n-1}} + \sum_{\mathbf{j}\in[-p+r+1,p-r-1]\times\{0\}^{n-1}} + \sum_{\mathbf{j}\in[p-r,p+r]\times\{0\}^{n-1}}\right](S_{\mathbf{j}}(\psi) - S_{\mathbf{j}}(U))$$

$$= \sum_{\mathbf{j}\in[-p-r,p+r]\times\{0\}^{n-1}} (S_{\mathbf{j}}(u) - S_{\mathbf{j}}(\psi))$$

$$(3.81) \quad + \left[\sum_{\mathbf{j}\in[-p-r,-p+r]\times\{0\}^{n-1}} + \sum_{\mathbf{j}\in[p-r,p+r]\times\{0\}^{n-1}}\right](S_{\mathbf{j}}(\psi) - S_{\mathbf{j}}(U)).$$

By Lemma 3.22, the first term on the right in (3.81) is $\leq 0$, while by (3.78), (3.79), the definition of $\psi$, and $S_{\mathbf{j}}$ is continuous on $\mathbb{R}^{\mathbb{Z}^n}$, the other terms on the right is $\leq \sigma/3$ in magnitude. On the other hand,

$$\sum_{\mathbf{j}\in[-p,p]\times\{0\}^{n-1}} (S_{\mathbf{j}}(u) - S_{\mathbf{j}}(U))$$

$$(3.82) \quad = J_{1;-p,p}(u) - J_{1;-p,p}(U)$$

$$\geq J_{1,-p,p}(u) - J_1(U) - \sigma/6$$

via (3.80). If $J_1(u) = \infty$, $J_{1,-p,p}(u) - J_1(U) - \sigma/6 \to \infty$ as $p \to \infty$. If $J_1(u) < \infty$, $J_{1,-p,p}(u) - J_1(U) - \sigma/6 \geq 2\sigma/3$ for large $p$. Both of the cases are contrary to (3.81). Thus we complete the proof of (3.75) and then Theorem 3.21. $\square$

## 4. Solutions heteroclinic in $\mathbf{i}_1$ and $\mathbf{i}_2$

In this section, we construct more complex heteroclinic solutions. We suppose $(*_0)$ holds and also $\mathcal{M}_1 = \mathcal{M}_1(v_0, w_0)$ has gaps, i.e.,

$(*_1)$ \qquad there are adjacent $v_1, w_1 \in \mathcal{M}_1(v_0, w_0)$ with $v_1 < w_1$.

Figure 3 illustrates these assumptions. We want to prove there is a solution lies between $v_1$ and $w_1$, which is heteroclinic in $\mathbf{i}_2$ from $v_1$ to $w_1$ as shown in Figure 4. The desired heteroclinic solution is periodic in $\mathbf{i}_3, \cdots, \mathbf{i}_n$. Since the proofs of theorems of this section are similar to Section 3, we will mainly state the results and omit the proofs.

Firstly, let $v, w \in \mathcal{M}_1$ with $v < w$. At this moment we do not require $v, w$ are adjacent in $\mathcal{M}_1$. Define

$$\hat{\Gamma}_2 = \hat{\Gamma}_2(v, w) = \{u \in \mathbb{R}^{\mathbb{Z}^2 \times (\mathbb{Z}/\{1\})^{n-2}} \mid v \leq u \leq w\}.$$

For $u \in \hat{\Gamma}_2$ and $l, i \in \mathbb{Z}$,

$$(4.1) \qquad \left\|\tau^2_{-i}\tau^1_{-l}u - v_0\right\|_{\mathbf{0}} \leq \|w - v_0\|_{\mathbf{T}_l} \to 0, \quad l \to -\infty$$

and similarly

$$(4.2) \qquad \left\|\tau^2_{-i}\tau^1_{-l}u - w_0\right\|_{\mathbf{0}} \to 0, \quad l \to \infty.$$

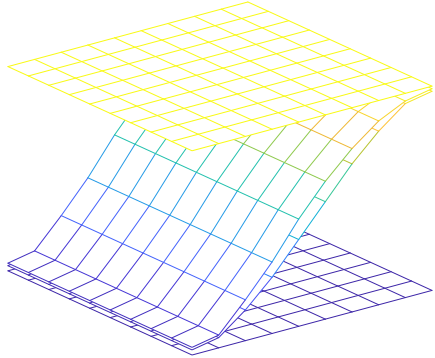

FIGURE 3. Assumptions of $(*_0)$ and $(*_1)$.

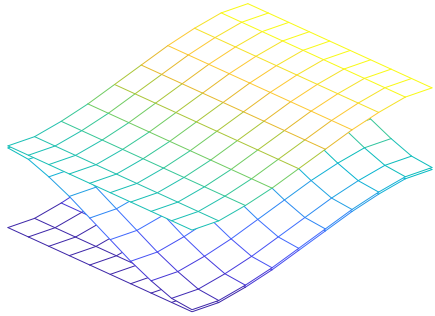

FIGURE 4. Heteroclinic solution in $\mathcal{M}_2(v_1, w_1)$.

Thus $\tau^2_{-i}u$ satisfies the asymptotic conditions in $\Gamma_1$ but $\tau^2_{-i}u \notin \Gamma_1$ because $\tau^2_{-i}u$ may be not periodic in $\mathbf{i}_2$. So $J_1(\tau^2_{-i}u)$ is not well-defined. We first extend the definition of $J_1$. As in Proposition 3.2, we have the following proposition.

**Proposition 4.1.** *For $u \in \hat{\Gamma}_2$, $J_{1;p,q}(u)$ is bounded from below and above independently of $u \in \hat{\Gamma}_2$ and $p, q$.*

**Proof:** For $u \in \hat{\Gamma}_2$, define $J_{1,i}(u), J_{1;p,q}(u)$ as before. Thus

$$J_{1;p,q}(u)$$

$$= \sum_{i=p}^{q} J_{1,i}(u)$$

(4.3)
$$= \sum_{i=p}^{q} S_{\mathbf{T}_i}(u) - (q + 1 - p)c_0$$

$$= \sum_{i=p}^{q} (S_{\mathbf{T}_i}(u) - S_{\mathbf{T}_i}(v)) + J_{1;p,q}(v).$$

As $p \to -\infty$, $q \to \infty$, $J_{1;p,q}(v) \to J_1(v) = c_1$. By Lemma 2.11,

$$|\sum_{i=p}^{q}(S_{\mathbf{T}_i}(u) - S_{\mathbf{T}_i}(v))|$$

$$\leq L \sum_{\mathbf{k} \in \overline{E_0 \cap \{p \leq \mathbf{i}_1 \leq q\}}} |u(\mathbf{k}) - v(\mathbf{k})|$$

(4.4)
$$\leq L \sum_{\mathbf{k} \in \overline{E_0 \cap \{p \leq \mathbf{i}_1 \leq q\}}} (w(\mathbf{k}) - v(\mathbf{k}))$$

$$\leq L(\#B_0^r) \sum_{\mathbf{k} \in E_0} (w(\mathbf{k}) - v(\mathbf{k})),$$

where $E_i := \mathbb{Z} \times \{i\} \times \{0\}^{n-2}$. Since $v, w \in \mathcal{M}_1$, $w < \tau_{-j}^1 v$ for some smallest $j > 0$. (Note we do not assume that $v < w$ are adjacent in $\mathcal{M}_1$.) Therefore

(4.5) $$\sum_{\mathbf{j} \in E_0} (w(\mathbf{j}) - v(\mathbf{j})) \leq \sum_{\mathbf{j} \in E_0} (\tau_{-j}^1 v(\mathbf{j}) - v(\mathbf{j})) \leq j(w_0 - v_0)(\mathbf{0}) \leq j.$$

This proves our proposition. $\square$

By (4.4), (4.5) and Cauchy criterion, we have

$$\lim_{\substack{p \to -\infty \\ q \to \infty}} J_{1;p,q}(u)$$

exists for $u \in \hat{\Gamma}_2$. We define

$$J_1(u) = \lim_{\substack{p \to -\infty \\ q \to \infty}} J_{1;p,q}(u),$$

thus by (4.3)

(4.6) $$J_1(u) = c_1 + \sum_{\mathbf{j} \in E_0} [S_{\mathbf{j}}(u) - S_{\mathbf{j}}(v)].$$

To construct solutions heteroclinic in $\mathbf{i}_1$ and $\mathbf{i}_2$, we need another renormalized functional $J_2(u)$. For $u \in \hat{\Gamma}_2$ and $i \in \mathbb{Z}$, set

$$J_{2,i}(u) = J_1(\tau_{-i}^2 u) - c_1,$$

$$J_{2;p,q}(u) = \sum_{i=p}^{q} J_{2,i}(u),$$

and
$$J_2(u) = \liminf_{\substack{p \to -\infty \\ q \to \infty}} J_{2;p,q}(u). \tag{4.7}$$

An analogue version of Proposition 3.2 for $J_{2;p,q}$ is ready, which ensures $J_2$ is well-defined.

**Proposition 4.2.** *Suppose $u \in \hat{\Gamma}_2(v,w)$ and $p,q \in \mathbb{Z}$. Then there is a constant $K_2 = K_2(v,w) \geq 0$ such that*
$$J_{2;p,q}(u) \geq -K_2.$$

**Proof:** By (4.4)-(4.5),
$$|J_{2,i}(u)| \leq L(\#B_\mathbf{0}^r) \sum_{\mathbf{k} \in E_i} (w(\mathbf{k}) - v(\mathbf{k})) = L(\#B_\mathbf{0}^r) \sum_{\mathbf{k} \in E_0} (w(\mathbf{k}) - v(\mathbf{k})) \leq L(\#B_\mathbf{0}^r)j =: M_2. \tag{4.8}$$

This proves the proposition for $q = p, p+1, \cdots, p+2r+2$ with any $K_2 \geq (2r+3)M_2$. Thus suppose $q > p+2r+2$ and define $\chi$ as in (3.5) with $\mathbf{i}_1$ replaced by $\mathbf{i}_2$. By Proposition 3.20, $J_{2;p,q}(\chi) \geq 0$. Continuing as in (3.6)-(3.7) (using $|u-v| \leq w(\mathbf{0}) - v(\mathbf{0})$) yields Proposition 4.2. □

The next lemma is similar to Lemma 3.3.

**Lemma 4.3.** *If $u \in \hat{\Gamma}_2$, $p,q \in \mathbb{Z}$ with $p \leq q$, then*
$$J_{2;p,q}(u) \leq J_2(u) + 2K_2. \tag{4.9}$$

It is useful to show that $J_{2,i}$ is continuous.

**Lemma 4.4.** *Suppose $\mathcal{Y} \subset \hat{\Gamma}_2$. Then $J_{2,i}$ is continuous (with respect to convergence in $\|\cdot\|_{E_{i-r} \cup \cdots \cup E_i \cup \cdots E_{i+r}}$) on $\mathcal{Y}$.*

**Proof:** For $u \in \mathcal{Y} \subset \hat{\Gamma}_2$, by (4.4) and (4.6) we have $J_{2,i}(u) < \infty$ for any $i \in \mathbb{Z}$. Let $(u_k) \subset \mathcal{Y}$, $u \in \mathcal{Y}$, and $\|u_k - u\|_{E_{i-r} \cup \cdots \cup E_i \cup \cdots E_{i+r}} \to 0$. By Lemma 2.11,
$$|J_{1;p,q}(\tau_{-i}^2 u_k) - J_{1;p,q}(\tau_{-i}^2 u)| \leq L \sum_{\mathbf{k} \in \bar{E}_i} |u_k(\mathbf{k}) - u(\mathbf{k})|.$$

Letting $p \to -\infty, q \to \infty$,
$$|J_{2,i}(u_k) - J_{2,i}(u)| \leq L \sum_{\mathbf{k} \in \bar{E}_i} |u_k(\mathbf{k}) - u(\mathbf{k})| \to 0 \tag{4.10}$$

as $k \to \infty$. □

Now similar to $\Gamma_1(v,w)$, we introduce
$$\Gamma_2 := \Gamma_2(v,w) := \{u \in \hat{\Gamma}_2 \mid \|u-v\|_{E_i} \to 0, i \to -\infty, \text{ and } \|u-w\|_{E_i} \to 0, i \to \infty\}.$$
As in Section 3, we have

**Proposition 4.5.** *For $u \in \Gamma_2$,*
$$\|\tau_{-i}^2 u - v\|_{E_0} \to 0, i \to -\infty, \tag{4.11}$$
$$\|\tau_{-i}^2 u - w\|_{E_0} \to 0, i \to \infty. \tag{4.12}$$
$$J_{2,i}(u) \to 0, |i| \to \infty. \tag{4.13}$$

*If $u \in \Gamma_2$ and $J_2(u) < \infty$, then*
$$J_2(u) = \lim_{\substack{p \to -\infty \\ q \to \infty}} J_{2;p,q}(u). \tag{4.14}$$

**Proof:** (4.11)-(4.12) follow from the definition of $\Gamma_2$ and (4.13) follows from (4.11)-(4.12) and Lemma 4.4. (4.14) is proved exactly as in the proof of (3.11) of Proposition 3.4. □

The next result for $J_2$ corresponds to Proposition 3.7.

**Proposition 4.6.** *Let* $\mathcal{Y} \subset \hat{\Gamma}_2(v, w)$. *Define*

$$(4.15) \qquad c(\mathcal{Y}) = \inf_{u \in \mathcal{Y}} J_2(u).$$

*If* $(u_k)$ *is a minimizing sequence for* (4.15), *then there is a* $U \in \hat{\Gamma}_2$ *such that along a subsequence,* $u_k \to U$ *in* $\mathbb{R}^{\mathbb{Z}^n}$ *and in* $\|\cdot\|_{E_j}$ *for any fixed* $j \in \mathbb{Z}$. *If* $c(\mathcal{Y}) < \infty$, *then*

$$(4.16) \qquad -K_2 \leq J_2(U) \leq c(\mathcal{Y}) + 1 + 2K_2,$$

*with* $K_2$ *as in Proposition 4.2.*

**Proof:** The existence of $U$ and $u_k \to U$ (may be up to a subsequence, which we still denote it by $u_k$) pointwise follow from Proposition 3.7. Note that

$$\sum_{\mathbf{j} \in E_i} |u_k(\mathbf{j}) - U(\mathbf{j})| \leq \sum_{\mathbf{j} \in E_i} |w(\mathbf{j}) - v(\mathbf{j})| \leq \bar{j},$$

for some $\bar{j}$ as in (4.5). Hence for $m = m(\epsilon) > 0$ sufficiently large, we have

$$\sum_{\mathbf{j} \in E_i \cap \{\|\mathbf{j}\| \geq m\}} |u_k(\mathbf{j}) - U(\mathbf{j})| \leq \epsilon.$$

Using $u_k \to U$ pointwise, we obtain

$$\sum_{\mathbf{j} \in E_i \cap \{\|\mathbf{j}\| \leq m\}} |u_k(\mathbf{j}) - U(\mathbf{j})| \to 0$$

as $k \to \infty$. Now Lemma 4.4 implies $J_{2;p,q}$ is continuous and thus (4.16) can be proved as (3.25). Then Proposition 4.6 follows. □

The following result is similar to Proposition 3.8. As before we need a notation. Fix $\mathbf{i} \in \mathbb{Z}^n$, Define $\delta_{\mathbf{i}}(\mathbf{j}) = 1$ for $\mathbf{j}_1 = \mathbf{i}_1$, and $\mathbf{j}_2 = \mathbf{i}_2$; $\delta_{\mathbf{i}}(\mathbf{j}) = 0$ otherwise.

**Proposition 4.7.** *Suppose there is a minimizing sequence* $(u_k)$ *for* (4.15) *such that for some* $\mathbf{i} \in \mathbb{Z}^n$, *the function* $\delta_{\mathbf{i}}$, *and* $t_0 > 0$, $c(\mathcal{Y}) \leq J_2(u_k + t\delta_{\mathbf{i}}) + \epsilon_k$ *for all* $|t| \leq t_0$, *where* $\epsilon_k \to 0$ *as* $k \to \infty$. *Then the limit* $U$ *of* $u_k$ *satisfies* (2.1) *at* $\mathbf{i}$. *Moreover* $U$ *satisfies* (2.1) *at any* $\mathbf{j}$ *with* $\mathbf{j}_1 = \mathbf{i}_1$ *and* $\mathbf{j}_2 = \mathbf{i}_2$.

**Proof:** As in Proposition 3.8 without any essential changes. □

The next theorem corresponding to Theorem 3.11. Before stating the theorem, we give some notations. For $v \in \mathcal{M}_1(v_0, w_0)$, set

$$\Gamma_2(v) = \{u \in \hat{\Gamma}_2(\tau_1^1 v, \tau_{-1}^1 v) \mid \|\tau_{-i}^2 u - v\|_{E_i} \to 0 \text{ as } |i| \to \infty\}.$$

Define

$$(4.17) \qquad c_2(v) = \inf_{u \in \Gamma_2(v)} J_2(u)$$

and let

$$\mathcal{M}_2(v) = \{u \in \Gamma_2(v) \mid J_2(u) = c_2(v)\}.$$

Then we obtain:

**Theorem 4.8.** *If* $s \in C^2(\mathbb{R}^{B_0^r}, \mathbb{R})$ *satisfies* (S1)-(S3) *and* $(*_0)$ *holds, then* $c_2(v) = 0$ *and* $\mathcal{M}_2(v) = \{v\}$.

**Proof:** The proof follows from the arguments of Theorem 3.11 with slight modifications.
□

Now we can state the second main theorem of this paper. Set

(4.18) $$c_2 = c_2(v_1, w_1) = \inf_{u \in \Gamma_2(v_1, w_1)} J_2(u),$$

and

$$\mathcal{M}_2 = \mathcal{M}_2(v_1, w_1) = \{u \in \Gamma_2(v_1, w_1) \mid J_2(u) = c_2\}.$$

**Theorem 4.9.** *If $s \in C^2(\mathbb{R}^{B_0^r}, \mathbb{R})$ satisfies (S1)-(S3) and $(*_0)$, $(*_1)$ holds, then there is a solution $U_2 \in \mathcal{M}_2$ of (2.1). Moreover, $\mathcal{M}_2$ is an ordered set and the elements of $\mathcal{M}_2$ are solutions of (2.1), and any $U \in \mathcal{M}_2$ is strictly 1-monotone in $\mathbf{i}_1, \mathbf{i}_2$.*

**Proof:** The proof is similar to the proof of Theorem 3.13 and we just point out the necessary modifications. The first one is (3.33) is replaced by

(4.19) $$u_k((0, i, 0, \cdots, 0)) \leq \frac{1}{2}(v_1 + w_1)(\mathbf{0}) \leq u_k(\mathbf{0})$$

for all $-i, k \in \mathbb{N}$. This will give a corresponding version of (3.35). (3.34) is replaced by

$$-K_2 \leq J_2(U_2) < c_2 + 1 + 2K_2.$$

Now proceeding as in the proof of Theorem 3.13, we should prove (A)-(D). For (A), the difference is (3.36) should be replaced by

$$\tau_2^1 w_1 \leq \tau_1^1 v_1 \leq u_k + t\delta_{\mathbf{i}} \leq \tau_{-1}^1 w_1$$

for $|t| \leq 1$. For (B), (3.41) is replaced by $U_2$ is 1-monotone in $\mathbf{i}_2$ and $U_2 \in \hat{\Gamma}_2(v_1, w_1) \setminus \{v_1, w_1\}$. Although Corollary 3.6 can not be used, a modified version, which is easy, would be enough to show that $U_2 \in \Gamma_2(v_1, w_1)$. (C) and (D) can be proved as in Theorem 3.13 with one exception:

(4.20) $$U < \tau_{-1}^1 U$$

for $U \in \mathcal{M}_2$.

To prove (4.20), define $\Phi = \max(U, \tau_{-1}^1 U)$ and $\Psi = \min(U, \tau_{-1}^1 U)$. We claim that

(4.21) $$\Phi \in \Gamma_2(\tau_{-1}^1 v_1, \tau_{-1}^1 w_1)$$

and

(4.22) $$\Psi \in \Gamma_2(v_1, w_1).$$

Suppose (4.21)-(4.22) for the moment. Thus we have

(4.23) $$J_2(\Phi) + J_2(\Psi) \leq J_2(U) + J_2(\tau_{-1}^1 U) = c_2(v_1, w_1) + c_2(\tau_{-1}^1 v_1, \tau_{-1}^1 w_1).$$

and then $J_2(\Phi) = c_2(\tau_{-1}^1 v_1, \tau_{-1}^1 w_1)$ and $J_2(\Psi) = c_2(v_1, w_1)$. Since the elements of $\mathcal{M}_2$ are solutions of (2.1), $\Phi$ and $\Psi$ are solutions of (2.1) with $\Phi \geq \Psi$. As earlier, we obtain

$$(i) \ U = \tau_{-1}^1 U, \quad \text{or} \quad (ii) \ U > \tau_{-1}^1 U, \quad \text{or} \quad (iii) \ U < \tau_{-1}^1 U.$$

But (i) contradicts $v_1 < U < w_1$ and (ii) leads to a contradiction:

$$w_0 > w_1 > U \geq \lim_{j \to \infty} \tau_{-j}^1 U \geq \lim_{j \to \infty} \tau_{-j}^1 v_1 = w_0.$$

Thus (iii) holds, which is (4.20).

What is left is to prove (4.21) - (4.22). Since the proofs of (4.21) and (4.22) being same, we only check (4.22). Since $v_1 < U$ and $v_1 < \tau_{-1}^1 v_1 < \tau_{-1}^1 U$, we have

$$v_1 < \Psi \leq U < w_1.$$

Therefore $\Psi \in \hat{\Gamma}_2$ and
$$\|\Psi - v_1\|_{E_i} \leq \|U - v_1\|_{E_i} \to 0, \quad i \to -\infty.$$

Next note that

(4.24) $$\sum_{E_i} |\Psi - w_1| = \sum_{E_i \cap \{\mathbf{i}_1 > r_1\}} |\Psi - w_1| + \sum_{E_i \cap \{\mathbf{i}_1 \leq r_1\}} |\Psi - w_1|.$$

Since $\left\|\tau^1_{-1} v_1 - v_1\right\|_{E_0} \leq \|w_0 - v_0\|_{\mathbf{0}} \leq 1$,

(4.25) $$\sum_{E_i \cap \{\mathbf{i}_1 > r_1\}} |\Psi - w_1| \leq \sum_{E_0 \cap \{\mathbf{i}_1 > r_1\}} (\tau^1_{-1} v_1 - v_1).$$

By Cauchy criterion, the right-hand side of (4.25) converges to 0 as $r_1 \to \infty$. Since $\tau^2_{-i} U \to w_1$ and $\tau^1_{-1}\tau^2_{-i} U \to \tau^1_{-1} w_1 > w_1$ as $i \to \infty$, convergence being in $\mathbb{R}^{E_0}$,

(4.26) $$\sum_{E_i \cap \{|\mathbf{i}_1| \leq r_1\}} |\Psi - w_1| \to 0, \quad i \to \infty.$$

Combining (4.24) – (4.26) gives
$$\|\Psi - w_1\|_{E_i} \to 0, \quad i \to \infty,$$

and the proof of Theorem 4.9 is complete. $\square$

Next as in Remark 3.15 we have:

**Remark 4.10.** *Suppose $s \in C^2(\mathbb{R}^{B^r_\mathbf{0}}, \mathbb{R})$ satisfy (S1) – (S3) and $(*_0)$ holds. Assume $v, w \in \mathcal{M}_1(v_0, w_0)$ with $v < w$. Then $\mathcal{M}_2(v, w) \neq \emptyset$ if and only if $v$ and $w$ are adjacent members of $\mathcal{M}_1(v_0, w_0)$. For the proof we refer to [22, Theorem 4.50].*

**Remark 4.11.** *Proposition 3.20 can be carried over to the current setting. The proof is simple and we omit it.*

Since the existence of elements of $\mathcal{M}_2$ depends on the gap conditions $(*_0)$ and $(*_1)$, we next explore the gap conditions $(*_0)$ and $(*_1)$. We have the following proposition which corresponds to Proposition 3.16.

**Proposition 4.12.** *Assume that $s \in C^2(\mathbb{R}^{B^r_\mathbf{0}}, \mathbb{R})$ satisfies (S1) – (S3) and $(*_0)$ and $(*_1)$ hold. Then there is an $\epsilon > 0$ such that for $\bar{s} \in C^2(\mathbb{R}^{B^r_\mathbf{0}}, \mathbb{R})$ satisfies (S1) – (S3), if*

(4.27) $$\|s - \bar{s}\|_{L^\infty(\mathbb{R}^{\#B^r_\mathbf{0}})} + \sum_{\|\mathbf{i}\| \leq r} \|\partial_\mathbf{i} s - \partial_\mathbf{i} \bar{s}\|_{L^\infty(\mathbb{R}^{\#B^r_\mathbf{0}})} \leq \epsilon$$

*is satisfied, $(*_0)$ and $(*_1)$ holds for $\bar{s}$. Moreover, suppose $v_1, w_1$ are a gap pair for $s$ for $(*_1)$ and*
$$\alpha_1 = v_1(\mathbf{0}); \qquad \beta_1 = w_1(\mathbf{0}).$$

*Then for any $\delta \in (0, \frac{\beta_1 - \alpha_1}{2})$, there is an $\epsilon_2 = \epsilon_2(s, \delta) > 0$ such that (3.58) holds with $\epsilon = \epsilon_2$ and*

(4.28) $$u(\mathbf{0}) \notin (\alpha_1 + \delta, \beta_1 - \delta)$$

*for all $u \in \mathcal{M}_1(v_0(\bar{s}), w_0(\bar{s}))$.*

**Proof:** By Proposition 3.16, for sufficiently small $\epsilon > 0$ $(*_0)$ holds for $\bar{s}$. Thus to complete the proof of Proposition 4.12, it suffices to verify (4.28). If (4.28) is false, there exist a $\delta \in (0, (\beta_1 - \alpha_1)/2)$ and a sequence $(s_k)$ satisfying (S1)-(S3) and (4.27) with $u_k \in \mathcal{M}_1(v_0(s_k), w_0(s_k))$ such that

$$u_k(\mathbf{0}) \in (\alpha_1 + \delta, \beta_1 - \delta). \tag{4.29}$$

Since $v_0(s_k) \leq u_k \leq w_0(s_k)$ and by Remark 3.17, $v_0(s_k), w_0(s_k)$ are near $v_0(s), w_0(s)$, it follows that $(u_k) \subset \hat{\Gamma}(v_0(s) - 1, w_0(s) + 1)$ for large $k$. By Proposition 3.7, (2.1) and (4.27) there is a solution $u$ of (2.1) for $s$ with

$$u(\mathbf{0}) \in [\alpha_1 + \delta, \beta_1 - \delta]. \tag{4.30}$$

The minimality of $u_k$ implies $u$ is minimal. Noting $\tau_{-1}^1 u_k > u_k$, we have

$$\tau_{-1}^1 u \geq u. \tag{4.31}$$

By Lemma 2.6, either $\tau_{-1}^1 u > u$ or $\tau_{-1}^1 u = u$. In both cases, $u$ is Birkhoff. Thus by (2) of Theorem 3.21, $u \in \mathcal{M}_0$ or $u \in \mathcal{M}_1(v_0(s), w_0(s))$. But this contradicts (4.30). Hence (4.28) holds and our proof is complete. $\square$

**Remark 4.13.** *Combining Propositions 3.16 and 4.12, if $\epsilon$ is small enough in (4.27), there will be a unique gap pair $v_1(\bar{s}), w_1(\bar{s})$ near $v_1(s), w_1(s)$.*

When $(*_1)$ does not hold, we can perturb $s$ to obtain $(*_1)$ again, as in Proposition 3.18.

**Theorem 4.14.** *Suppose $s \in C^2(\mathbb{R}^{B_0^r}, \mathbb{R})$ satisfy (S1) – (S3). Then for any $\epsilon > 0$, there is an $\bar{s} \in C^2(\mathbb{R}^{B_0^r}, \mathbb{R})$ satisfying (S1) – (S3), (4.27) and $(*_0), (*_1)$ with*
*(1) $\mathcal{M}_0(\bar{s}) = \{v + j \mid j \in \mathbb{Z}\}$ for some prescribed $v \in \mathcal{M}_0(s)$.*
*(2) $\mathcal{M}_1(v, v + 1, \bar{s}) = \{\tau_{-k}^1 U \mid k \in \mathbb{Z}\}$ for some $U \in \mathcal{M}_1(v, v+1, \bar{s})$, where $\mathcal{M}_1(v, v+1, \bar{s})$ is the set of minimizers given by Theorem 3.13.*

**Proof:** For prescribed $v \in \mathcal{M}_0(s)$, define

$$\bar{s}_1(u) = \sin^2 \pi[u(\mathbf{0}) - v(\mathbf{0})] + \frac{1}{2} \sum_{\mathbf{k}: \|\mathbf{k}\| \leq r} (u(\mathbf{k}) - u(\mathbf{0}))^2.$$

As in the proof of Proposition 3.18, we have $s + \delta_1 \bar{s}_1$ satisfies (S1) – (S3) and (4.27), (1) of Theorem 4.14, $(*_0)$ hold. Depending on whether $(*_1)$ holds for $s + \delta_1 \bar{s}_1$ or not, our arguments are divided into two parts.

(A). Assume that $(*_1)$ holds for $s + \delta_1 \bar{s}_1$. If (2) of Theorem 4.14 holds for $\bar{s} = s + \delta_1 \bar{s}_1$, we are through. Thus suppose (2) of Theorem 4.14 dose not hold. Fix any $U_1 \in \mathcal{M}_1(v, v+1, s + \delta_1 \bar{s}_1)$. We define a functional $\bar{s}_2$ as follows. For $u \in \mathbb{R}^{\mathbb{Z}^n}$, if

$$u(\mathbf{0}) \in \{v(\mathbf{0}) + k \mid k \in \mathbb{Z}\} \cup \{\tau_{-j}^1 U_1(\mathbf{0}) + k \mid j \in \mathbb{Z}, k \in \mathbb{Z}\},$$

then define $\bar{s}_2(u) = 0$; if not, then there exist $j \in \mathbb{Z}$ and $k \in \mathbb{Z}$ such that

$$u(\mathbf{0}) \in (\tau_{-j}^1 U_1(\mathbf{0}) + k, \tau_{-j-1}^1 U_1(\mathbf{0}) + k).$$

Define

$$\bar{s}_2(u) = |u(\mathbf{0}) - \tau_{-j}^1 U_1(\mathbf{0}) - k|^4 |u(\mathbf{0}) - \tau_{-j-1}^1 U_1(\mathbf{0}) - k|^4.$$

Noticing $\mathcal{M}_1(v, v+1, s + \delta_1 \bar{s}_1)$ is an ordered set, we have $\bar{s}_2(u) > 0$ on

$$\mathcal{M}_1(v, v+1, s + \delta_1 \bar{s}_1) \setminus \{\tau_{-j}^1 U_1 \mid j \in \mathbb{Z}\}.$$

Set $\bar{s} = s + \delta_1 \bar{s}_1 + \delta_2 \bar{s}_2$.

Now for $\delta_1, \delta_2 > 0$ small enough, we claim $\bar{s}$ is a desired functional. Indeed, (S1)-(S3) are easy to verify, so we only check (4.27), $(*_0)$, $(*_1)$ and (1)-(2) of Theorem 4.14. Since $\delta_1, \delta_2 > 0$ small enough, (4.27) holds. Note $J_0^{\bar{s}}(v) = c_0(s)$ and if $u \in \Gamma_0 \setminus \{v + j \,|\, j \in \mathbb{Z}\}$,
$$J_0^{\bar{s}}(u) > J_0^s(u) \geq c_0(s).$$
So $c_0(s) = c_0(\bar{s})$. This proves (1) of Theorem 4.14 and thus $(*_0)$. Similarly if $u \in \Gamma_1(v, v+1)$,
$$J_1^{\bar{s}}(u) = J_1^{s+\delta_1 \bar{s}_1}(u) + \delta_2 \sum_{\mathbf{j} \in \mathbb{Z} \times \{\mathbf{0}\}^{n-1}} \bar{s}_2(\tau_{\mathbf{j}_1}^1 u|_{B_{\mathbf{0}}^r}) \geq c_1(s + \delta_1 \bar{s}_1).$$
So
$$(4.32) \qquad c_1(s + \delta_1 \bar{s}_1) = c_1(\bar{s}).$$
Moreover, if $u \in \mathcal{M}_1(\bar{s}) \setminus \{\tau_{-j}^1 U_1 \,|\, j \in \mathbb{Z}\}$, then $J_1^{\bar{s}}(u) = c_1(\bar{s}) > c_1(s + \delta_1 \bar{s}_1)$, contrary to (4.32). This proves (2) of Theorem 4.14 and $(*_1)$. Our claim is proved.

(B). Assume that $(*_1)$ does not hold for $s + \delta_1 \bar{s}_1$. We can still define $\bar{s}_2$ and prove that $\bar{s} := s + \delta_1 \bar{s}_1 + \delta_2 \bar{s}_2$ is the desired functional. This completes the proof of Theorem 4.14. □

The rest of this section is devoted to the relation between the elements in $\mathcal{M}_2$ and minimal and Birkhoff solutions. First we have:

**Theorem 4.15.** *Assume that $s \in C^2(\mathbb{R}^{B_0^r}, \mathbb{R})$ satisfies (S1) - (S3) and $(*_0)$, $(*_1)$ hold. Then every element in $\mathcal{M}_2(v_1, w_1) \cup \mathcal{M}_2(w_1, v_1)$ is minimal and Birkhoff.*

The proof depends on Remark 4.11 and (1) of Theorem 3.21, which is easy and we omit it.

One may expect to extend (2) of Theorem 3.21 directly, but it is not true. In fact, since the dimension is higher there are more possibilities for $u$. But we have:

**Proposition 4.16.** *Assume that $s \in C^2(\mathbb{R}^{B_0^r}, \mathbb{R})$ satisfies (S1) - (S3) and set $(*_0)$, $(*_1)$ hold. If $U \in \Gamma_2(v_1, w_1)$ is minimal and Birkhoff, then $U \in \mathcal{M}_2(v_1, w_1)$.*

**Proof:** $U$ is a solution of (2.1) since it is minimal. Since $U \in \Gamma_2(v_1, w_1)$, we need only to check
$$(4.33) \qquad J_2(U) = c_2.$$
To verify (4.33), we use an argument analogous to the proof of Theorem 3.21. First note that as $j \to \infty$,
$$(4.34) \qquad \|U - w_1\|_{E_j}, \quad \|U - v_1\|_{E_{-j}} \to 0.$$
If (4.33) is not true, we have
$$(4.35) \qquad J_2(U) > c_2.$$
Choose $\psi \in \Gamma_2(v_1, w_1)$ such that for some $\sigma > 0$,
$$(4.36) \qquad c_2 \leq J_2(\psi) < J_2(\psi) + \sigma < J_2(U).$$
By (4.34) and (4.11)-(4.12), for any $\kappa > 0$, there is a $q = q(\kappa) \in \mathbb{N}$ such that for $\phi \in \{U, \psi\}$,
$$(4.37) \qquad \begin{aligned} \|\phi - v_1\|_{E_i} &\leq \kappa, \quad i \leq -q, \\ \|\phi - w_1\|_{E_i} &\leq \kappa, \quad i \geq q. \end{aligned}$$
For $i \in \mathbb{Z}$ and $k = -r, \cdots, 0, \cdots, r$, set
$$G_i^k = \begin{cases} U, & \mathbf{i}_2 \leq i + k, \\ \psi, & \mathbf{i}_2 \geq i + k + 1, \end{cases}$$

$$H_i^k = \begin{cases} U, & \mathbf{i}_2 \geq i+k+1, \\ \psi, & \mathbf{i}_2 \leq i+k. \end{cases}$$

Thus by (4.37) and Lemma 4.4, for $\kappa = \kappa(\sigma)$ small enough and $\phi \in \{U, \psi\} \cup \{G_i^k, H_i^k\}_{k=-r}^r$,

$$(4.38) \qquad |J_{2,i}(\phi)| \leq \frac{\sigma}{16(2r+1)}$$

for $|i| > q(\kappa) + 2r$. For $p \in \mathbb{N}$ and $p \geq q(\kappa) + r + 1$ large enough,

$$(4.39) \qquad J_{2,-p,p}(\psi) \leq J_2(\psi) + \frac{\sigma}{6}.$$

Setting

$$\Psi = \begin{cases} \psi, & -p+1 \leq \mathbf{i}_2 \leq p-1, \\ U, & \text{otherwise}, \end{cases}$$

we have

$$\sum_{\mathbf{j} \in \mathbb{Z} \times [-p-r, p+r] \times \{0\}^{n-2}} [S_{\mathbf{j}}(U) - S_{\mathbf{j}}(\psi)]$$

$$= \sum_{\mathbf{j} \in \mathbb{Z} \times [-p-r, p+r] \times \{0\}^{n-2}} [S_{\mathbf{j}}(U) - S_{\mathbf{j}}(\Psi)] + \sum_{\mathbf{j} \in \mathbb{Z} \times [-p-r, p+r] \times \{0\}^{n-2}} [S_{\mathbf{j}}(\Psi) - S_{\mathbf{j}}(\psi)]$$

$$= \sum_{\mathbf{j} \in \mathbb{Z} \times [-p-r, p+r] \times \{0\}^{n-2}} [S_{\mathbf{j}}(U) - S_{\mathbf{j}}(\Psi)]$$

$$+ \Big( \sum_{\mathbb{Z} \times [-p-r, -p+r] \times \{0\}^{n-2}} + \sum_{\mathbb{Z} \times [-p+r+1, p-r-1] \times \{0\}^{n-2}} + \sum_{\mathbb{Z} \times [p-r, p+r] \times \{0\}^{n-2}} \Big)[S_{\mathbf{j}}(\Psi) - S_{\mathbf{j}}(\psi)]$$

$$= \sum_{\mathbf{j} \in \mathbb{Z} \times [-p-r, p+r] \times \{0\}^{n-2}} [S_{\mathbf{j}}(U) - S_{\mathbf{j}}(\Psi)]$$

$$+ \Big( \sum_{\mathbb{Z} \times [-p-r, -p+r] \times \{0\}^{n-2}} + \sum_{\mathbb{Z} \times [p-r, p+r] \times \{0\}^{n-2}} \Big)[S_{\mathbf{j}}(\Psi) - S_{\mathbf{j}}(\psi)]$$

$$= \sum_{\mathbf{j} \in \mathbb{Z} \times [-p-r, p+r] \times \{0\}^{n-2}} [S_{\mathbf{j}}(U) - S_{\mathbf{j}}(\Psi)]$$

$$+ \sum_{i=-p-r}^{-p+r} [J_{2,i}(\Psi) - J_{2,i}(\psi)] + \sum_{i=p-r}^{p+r} [J_{2,i}(\Psi) - J_{2,i}(\psi)]$$

$$= \sum_{\mathbf{j} \in \mathbb{Z} \times [-p-r, p+r] \times \{0\}^{n-2}} [S_{\mathbf{j}}(U) - S_{\mathbf{j}}(\Psi)]$$

$$(4.40) \quad + \sum_{i=-p-r}^{-p+r} [J_{2,i}(G_i^{-p-i}) - J_{2,i}(\psi)] + \sum_{i=p-r}^{p+r} [J_{2,i}(H_i^{p-i}) - J_{2,i}(\psi)].$$

The first term on the right is $\leq 0$, which follows from Lemma 4.17 below. Lemma 4.17 is an analogue of Lemma 3.22. By (4.38), for large $p$, each of the remaining terms on the right is $\leq \frac{\sigma}{16(2r+1)}$ in magnitude. To estimate the left-hand side of (4.40), we write

$$\sum_{\mathbf{j} \in \mathbb{Z} \times [-p-r, p+r] \times \{0\}^{n-2}} (S_{\mathbf{j}}(U) - S_{\mathbf{j}}(\psi))$$

$$= \sum_{i=-p-r}^{p+r} \sum_{\mathbf{j} \in E_i} (S_{\mathbf{j}}(U) - S_{\mathbf{j}}(\psi))$$

(4.41)
$$= \sum_{i=-p-r}^{p+r} (J_1(\tau_{-i}^2 U) - J_1(\tau_{-i}^2 \psi))$$

$$= J_{2;-p-r,p+r}(U) - J_{2;-p-r,p+r}(\psi)$$

$$\geq J_{2;-p-r,p+r}(U) - J_2(\psi) - \frac{\sigma}{6}$$

via (4.39). Thus if $J_2(U) = \infty$, by (4.41) the left-hand side of (4.40)$\to \infty$ as $p \to \infty$, while if $J_2(U) < \infty$, by (4.36), the left-hand side of (4.40) exceeds $2\sigma/3$. Either case leads to a contradiction and (4.33) is valid, which completes the proof of Proposition 4.16. □

To complete the proof of Proposition 4.16, we need to prove

**Lemma 4.17.** *If $u \in \Gamma_2(v_1, w_1)$ is minimal, then for any $\phi \in \Gamma_2(v_1, w_1)$ with $\mathrm{supp}(\phi - u) \subset \mathbb{Z} \times [p, q] \times \mathbb{Z}^{n-2}$,*

(4.42)
$$\sum_{\mathbb{Z}^2 \times \{0\}^{n-2}} [S_{\mathbf{j}}(\phi) - S_{\mathbf{j}}(u)] \geq 0,$$

*where $p, q \in \mathbb{Z}$ and $p \leq q$.*

**Proof:** Let $\psi = \phi - u$ and define $\theta_l$ as in Lemma 3.22. Suppose $n = 2$ and $m \in \mathbb{N}$. Define $(\theta_m^1 \psi)(\mathbf{j}) = \theta_m(|\mathbf{j}_1|)\psi(\mathbf{j})$. Since $u$ is minimal,

$$0 \leq \sum_{\mathbb{Z}^2} [S_{\mathbf{j}}(u + \theta_m^1 \psi) - S_{\mathbf{j}}(u)]$$

(4.43)
$$= \sum_{[-m-r, m+r] \times [p-r, q+r]} [S_{\mathbf{j}}(u + \theta_m^1 \psi) - S_{\mathbf{j}}(u)]$$

$$= \sum_{[-m+r, m-r] \times [p-r, q+r]} [S_{\mathbf{j}}(u + \psi) - S_{\mathbf{j}}(u)] + \mathcal{R}_1(u, \psi) + \mathcal{R}_2(u, \psi),$$

where

$$\mathcal{R}_1(u, \psi) = \sum_{A_1} [S_{\mathbf{j}}(u + \theta_m^1 \psi) - S_{\mathbf{j}}(u)],$$

$$\mathcal{R}_2(u, \psi) = \sum_{A_2} [S_{\mathbf{j}}(u + \theta_m^1 \psi) - S_{\mathbf{j}}(u)],$$

and $A_1, A_2$ are the regions

$$A_1 = [-m-r, -m+r-1] \times [p-r, q+r],$$
$$A_2 = [m-r+1, m+r] \times [p-r, q+r].$$

As in (4.4)-(4.5), we have

$$\mathcal{R}_1(u, \psi), \quad \mathcal{R}_2(u, \psi) \to 0, \qquad \text{as } m \to \infty,$$

since they are bounded by tails of a convergence series that does not depend on $m$. Now letting $m \to \infty$ in (4.43) yields (4.42).

If $n > 2$, define $(\theta_m^1 \theta^2 \psi)(\mathbf{j}) = \theta_m(|\mathbf{j}_1|)\theta_l(|\mathbf{j}_3|) \cdots \theta_l(|\mathbf{j}_n|)\psi(\mathbf{j})$. As in (4.43),

$$
\begin{aligned}
0 &\leq \sum_{\mathbb{Z}^n}[S_{\mathbf{j}}(u + \theta_m^1 \theta^2 \psi) - S_{\mathbf{j}}(u)] \\
&= \sum_{[-m-r,m+r] \times [p-r,q+r] \times [-l-r,l+r]^{n-2}}[S_{\mathbf{j}}(u + \theta_m^1 \theta^2 \psi) - S_{\mathbf{j}}(u)] \\
&= (2l - 2r + 1)^{n-2} \sum_{[-m+r,m-r] \times [p-r,q+r] \times \{0\}^{n-2}}[S_{\mathbf{j}}(u + \psi) - S_{\mathbf{j}}(u)] \\
&\quad + \mathcal{R}_{1,l}(u, \psi) + \mathcal{R}_{2,l}(u, \psi),
\end{aligned}
\tag{4.44}
$$

where

$$
\mathcal{R}_{1,l}(u, \psi) = \sum_{A_{1,l}}[S_{\mathbf{j}}(u + \theta_m^1 \theta^2 \psi) - S_{\mathbf{j}}(u)],
$$

$$
\mathcal{R}_{2,l}(u, \psi) = \sum_{A_{2,l}}[S_{\mathbf{j}}(u + \theta_m^1 \theta^2 \psi) - S_{\mathbf{j}}(u)],
$$

and $A_{1,l}, A_{2,l}$ are the regions

$$
A_{1,l} = [-m-r, -m+r-1] \times [p-r, q+r] \times ([-l-r, l+r]^{n-2} \setminus [-l+r, l-r]^{n-2}),
$$

$$
A_{2,l} = [m-r+1, m+r] \times [p-r, q+r] \times ([-l-r, l+r]^{n-2} \setminus [-l+r, l-r]^{n-2}).
$$

Proceeding as in (3.67), first letting $l \to \infty$ and then $m \to \infty$ yields (4.42). □

The final proposition concerns the relations between heteroclinic solutions in $\mathcal{M}_2$ obtained for $s$ and $s_k$, which are perturbed potentials of $s$. Suppose $s, s_k \in C^2(\mathbb{R}^{B_0^r}, \mathbb{R})$ satisfy (S1)-(S3) and (4.27) holds for $\epsilon = \epsilon_k \to 0$ as $k \to \infty$. Assume that $(*_0), (*_1)$ hold for $s$. By Remark 4.13, $v_1(s_k), w_1(s_k)$ approach $v_1(s), w_1(s)$ as $k \to \infty$, respectively. Let $U_k \in \mathcal{M}_2(v_1(s_k), w_1(s_k))$ be given by Theorem 4.9. We have

**Proposition 4.18.** *There is a subsequence of $(U_k)$ converging to some point of*

$$\mathcal{M}_2(v_1(s), w_1(s)) \cup \{v_1(s), w_1(s)\}.$$

For the proof of Proposition 4.18, please see [22, Proposition 4.79].

## 5. Three generalizations

In this section, we generalize the results of Sections 3-4. In Section 5.1, we construct heteroclinic solutions in directions $\mathbf{i}_1, \cdots, \mathbf{i}_k$, $3 \leq k \leq n$. Section 5.2 deals with other linearly independent directions other than $\mathbf{e}_1, \cdots, \mathbf{e}_n$. We loose the assumption $\alpha = \mathbf{0}$ in Section 5.3 and extend our results to $\alpha \in \mathbb{Q}^n \setminus \{\mathbf{0}\}$.

**5.1. Higher dimensional heteroclinic solutions.** Recalling when $(*_0)$ holds, we construct two types of heteroclinic solutions, $\mathcal{M}_1(v_0, w_0)$ and $\mathcal{M}_1(w_0, v_0)$. Thus $(*_1)$ has two versions. Similarly for each version of $(*_1)$, there are two types of heteroclinic solutions, $\mathcal{M}_2(v_1, w_1)$ and $\mathcal{M}_2(w_1, v_1)$. By induction, for $1 \leq k \leq n$, we have $2^k$ versions of $(*_k)$. But in an obvious manner, it is enough to deal with one version of $(*_k)$. We treat the case of $U \in \mathcal{M}_k(v_{k-1}, w_{k-1})$ with $\tau_{-1}^i U > U$, $1 \leq i \leq k$.

We prove our theories by induction so suppose they hold for $l < n$. For $l + 1$, assume that

($*_l$)         there are adjacent $v_l, w_l \in \mathcal{M}_l(v_{l-1}, w_{l-1}) =: \mathcal{M}_l$ with $v_l < w_l$.

For $v, w \in \mathcal{M}_l$ with $v < w$ define

$$\hat{\Gamma}_{l+1} := \hat{\Gamma}_{l+1}(v, w) := \{u \in \mathbb{R}^{\mathbb{Z}^{l+1} \times (\mathbb{Z}/\{1\})^{n-(l+1)}} \mid v \leq u \leq w\}.$$

As before, for $u \in \hat{\Gamma}_{l+1}$ and $i \in \mathbb{Z}$, the functions $\tau^{l+1}_{-i} u$ have asymptotic limits in the directions $\mathbf{i}_j$, $1 \leq j \leq l$, but $J_l(\tau^{l+1}_{-i} u)$ is not well-defined. Setting $E^{l+1}_i = \mathbb{Z}^l \times \{i\} \times \{0\}^{n-(l+1)}$ and replacing $E_0$ of Section 4 by $E^{l+1}_0$ shows how $J_l$ extends to this setting and as in (4.6),

(5.1) $$J_l(u) = c_l + \sum_{\mathbf{j} \in E^{l+1}_0} [S_\mathbf{j}(u) - S_\mathbf{j}(v)].$$

Thus we can define $J_{l+1,i}(u)$ for $u \in \hat{\Gamma}_{l+1}$ via

$$J_{l+1,i}(u) := J_l(\tau^{l+1}_{-i} u) - c_l.$$

Following Section 4 we define $J_{l+1}$ as

$$J_{l+1}(u) := \liminf_{\substack{p \to -\infty \\ q \to \infty}} J_{l+1;p,q}(u) := \liminf_{\substack{p \to -\infty \\ q \to \infty}} \sum_{i=p}^{q} J_{l+1,i}(u).$$

Lemma 4.3 has an extension in this setting. Letting

$$\Gamma_{l+1} := \Gamma_{l+1}(v, w)$$
$$:= \{u \in \hat{\Gamma}_{l+1} \mid \|u - v\|_{E^{l+1}_i} \to 0, i \to -\infty; \|u - w\|_{E^{l+1}_i} \to 0, i \to \infty\}$$

leads to extensions of Lemma 4.4, Propositions 4.6 and 4.7, and Theorem 4.8. Set

(5.2) $$c_{l+1} := c_{l+1}(v_l, w_l) := \inf_{u \in \Gamma_{l+1}(v_l, w_l)} J_{l+1}(u),$$

and

$$\mathcal{M}_{l+1} := \mathcal{M}_{l+1}(v_l, w_l) := \{u \in \Gamma_{l+1}(v_l, w_l) \mid J_{l+1}(u) = c_{l+1}\}.$$

With these notations, we have:

**Theorem 5.1.** *If $s \in C^2(\mathbb{R}^{B^r_0}, \mathbb{R})$ satisfies (S1)-(S3) and $(*_i)$ holds for $i = 0, \cdots, l$, then there is a solution $U_{l+1} \in \mathcal{M}_{l+1}$ of (2.1). Moreover, $\mathcal{M}_{l+1}$ is an ordered set and the elements of $\mathcal{M}_{l+1}$ are solutions of (2.1), and any $U \in \mathcal{M}_{l+1}$ is strictly 1-monotone in $\mathbf{i}_1, \cdots, \mathbf{i}_{l+1}$.*

Remark 4.10 has an extension of the following form.

**Remark 5.2.** *Assume that $s \in C^2(\mathbb{R}^{B^r_0}, \mathbb{R})$ satisfies (S1)-(S3) and $(*_i)$ holds, $i = 0, \cdots, l-1$, and $v, w \in \mathcal{M}_l$ with $v < w$. Then $\mathcal{M}_{l+1}(v, w) \neq \emptyset$ if and only if $v$ and $w$ are adjacent members of $\mathcal{M}_l$.*

Proposition 4.12 and Theorem 4.14 are extended to the case of $(*_l)$. Theorem 4.15 is extended to:

**Theorem 5.3.** *Assume that $s \in C^2(\mathbb{R}^{B^r_0}, \mathbb{R})$ satisfies (S1)-(S3) and $(*_i)$ holds, $i = 0, \cdots, l-1$. If $u \in \mathcal{M}_l(v_{l-1}, w_{l-1})$, then $u$ is minimal and Birkhoff.*

**Proof:** The proof is as in Theorem 4.15. □

Proposition 4.18 also has an extension here. We do not list it since the generalization is easy.

**5.2. Other coordinate systems.** Suppose $\omega_i = \sum_{j=1}^{n} \alpha_{ij} \mathbf{e}_j$ with $\alpha_{ij} \in \mathbb{Z}, 1 \leq i, j \leq n$, and the vectors $\omega_i$ are linearly independent. Without loss of generality, we assume that the $\omega_i$ are orthogonal and for fixed $i$, the components $\alpha_{ij}$ of $\omega_i$ have no common factor. Replacing $\mathbf{e}_i$ by $\omega_i$, we construct heteroclinic solutions corresponding to $\omega_i$. Firstly let us consider periodic solutions, i.e., $u(\mathbf{i} + \omega_i) = u(\mathbf{i}), 1 \leq i \leq n$. To this end, set $\omega := (\omega_1, \cdots, \omega_n)$,

$$\mathcal{R} := \mathcal{R}(\omega) := \mathbb{Z}^n \cap \Big\{ \sum_{i=1}^{m} t_i \omega_i \,\big|\, 0 \leq t_i < 1,\, 1 \leq i \leq n \Big\}$$

and

(5.3) $$\Gamma_0(\omega) := \{ u \in \mathbb{R}^{\mathbb{Z}^n} \,|\, u(\mathbf{i} + \omega_i) = u(\mathbf{i}), 1 \leq i \leq n \}.$$

For $u \in \Gamma_0(\omega)$, let

(5.4) $$J_0^\omega(u) := \sum_{\mathbf{j} \in \mathcal{R}} S_{\mathbf{j}}(u)$$

and set

(5.5) $$c_0(\omega) := \inf_{u \in \Gamma_0(\omega)} J_0^\omega(u).$$

Following Theorem 2.12, we have a set $\mathcal{M}_0(\omega)$ of minimizers of this variational problem and $\mathcal{M}_0(\omega)$ is ordered. Moreover, the results in Sections 3 – 4 and Section 5.1 can be generalized in this setting with $\mathbf{e}_1, \cdots, \mathbf{e}_n$ replaced by $\omega$. However, they do not produce more new solutions as one may expect. For example:

**Lemma 5.4.** $\mathcal{M}_0(\omega) = \mathcal{M}_0(\mathbf{e}_1, \cdots, \mathbf{e}_n)$.

Comparing to the proof of [22, Lemma 5.9], there is nothing new in the proof of Lemma 5.4 so we omit it. When $(*_0)$ holds, as in Section 3, we can define $\Gamma_1(v_0, w_0; \omega)$, $c_1(v_0, w_0; \omega)$ and $\mathcal{M}_1(v_0, w_0; \omega)$. $\mathcal{M}_1(v_0, w_0; \omega)$ is an ordered set and it satisfies the following proposition.

**Proposition 5.5.** *Let $\omega = (\omega_1, \cdots, \omega_n)$ and $\hat{\omega} = (\hat{\omega}_1, \cdots, \hat{\omega}_n)$ be admissible sets of orthogonal vectors. Then*

$$\mathcal{M}_1(v_0, w_0; \omega) = \mathcal{M}_1(v_0, w_0; \hat{\omega}) \iff \omega_1 = \hat{\omega}_1.$$

We omit the proof since it is same to [22, Proposition 5.11]. By Proposition 5.5, $\mathcal{M}_1(v_0, w_0; \cdot)$ is determined by $\omega_1$, thus we denote $\mathcal{M}_1(v_0, w_0; \omega)$ by $\mathcal{M}_1(v_0, w_0; \omega_1)$. As one may expect, $\mathcal{M}_2(v_1, w_1; \omega)$ may only depend on $\omega_2$. This is the case, as the next proposition shows. Before stating the proposition, we need a remark on the notations. Since $u \in \mathcal{M}_1(v_0, w_0; \omega_1)$ only depends on $\omega_1$, we will denote $u$ by $u(\omega_1)$ and thus the associated $\mathcal{M}_2(v_1, w_1; \omega)$ by $\mathcal{M}_2(v_1(\omega_1), w_1(\omega_1); \omega)$.

**Proposition 5.6.** *Let $\omega = (\omega_1, \omega_2, \cdots, \omega_n)$ and $\hat{\omega} = (\omega_1, \hat{\omega}_2, \cdots, \hat{\omega}_n)$. Then*

$$\mathcal{M}_2(v_1(\omega_1), w_1(\omega_1); \omega) = \mathcal{M}_2(v_1(\omega_1), w_1(\omega_1); \hat{\omega}) \iff \omega_2 = \hat{\omega}_2.$$

**Proof:** The necessity follows as in the proof of Proposition 5.5. Now suppose $\omega_2 = \hat{\omega}_2$ and $u \in \mathcal{M}_2(v_1(\omega_1), w_1(\omega_1); \omega)$. As in the proof of sufficiency of [22, Proposition 5.11], we have $u(\mathbf{i} + \hat{\omega}_i) = u(\mathbf{i})$, $3 \leq i \leq n$, so $u \in \hat{\Gamma}_2(v_1(\omega_1), w_1(\omega_1))$. We claim that $u \in \Gamma_2(v_1(\omega_1), w_1(\omega_1); \hat{\omega})$, which needs to prove

(5.6) $$\|u - v_1(\omega_1)\|_{E_i^{\hat{\omega}}} \to 0, \, i \to -\infty,$$

(5.7) $$\|u - w_1(\omega_1)\|_{E_i^{\hat{\omega}}} \to 0, \, i \to \infty.$$

Here $E_i^{\hat{\omega}}$ is similar to $E_i$. Thus $E_i^{\hat{\omega}} = E_0^{\hat{\omega}} + i\omega_2$ and

$$E_0^{\hat{\omega}} = \left\{ t_1 \omega_1 + t_2 \omega_2 + \sum_{i=3}^{n} t_i \hat{\omega}_i \middle| t_1 \in \mathbb{R}, \, 0 \leq t_i < 1, 2 \leq i \leq n \right\} \cap \mathbb{Z}^n.$$

Note that

$$E_0^{\hat{\omega}}$$
$$= \left\{ t_1 \omega_1 + t_2 \omega_2 + \sum_{i,k=3}^{n} t_i q_{ik} \omega_k \middle| t_1 \in \mathbb{R}, \, 0 \leq t_i < 1, 2 \leq i \leq n \right\} \cap \mathbb{Z}^n$$
$$\subset \left\{ \sum_{i=1}^{n} t_i \omega_i \middle| t_1 \in \mathbb{R}, 0 \leq t_2 < 1, |t_i| \leq j, 3 \leq i \leq n \right\} \cap \mathbb{Z}^n$$
$$=: E^*$$

for some $j \in \mathbb{N}$. Hence

(5.8) $$\|u - v_1(\omega_1)\|_{E_i^{\hat{\omega}}} \leq \|u - v_1(\omega_1)\|_{E^* + i\omega_2},$$

and since $u \in \Gamma_2(v_1(\omega_1), w_1(\omega_1); \omega)$, we have $\|u - v_1(\omega_1)\|_{E^* + i\omega_2} \to 0$ as $i \to -\infty$. Thus (5.6) is satisfied. Similarly we obtain (5.7). Consequently, $u \in \Gamma_2(v_1(\omega_1), w_1(\omega_1); \hat{\omega})$. Since $u$ is also minimal and Birkhoff, by a variant of Proposition 4.16,

$$u \in \mathcal{M}_2(v_1(\omega_1), w_1(\omega_1); \hat{\omega})$$

and the proof of Proposition 5.6 is complete. □

More heteroclinic solutions of (2.1) as in Section 5.1 can also be constructed but we omit it here since the extension is easy.

**5.3. Generalizations to $\alpha \in \mathbb{Q}^n$.** The aim of this subsection is to replace the condition $\alpha = \mathbf{0}$ by $\alpha \in \mathbb{Q}^n$. The first step is to obtain $\mathcal{M}_0$. In [15], Moser constructed minimal solutions $u^*$ without self-intersections of rational rotation vector $\alpha \in \mathbb{Q}^n$ by minimization method. He translated searching such solutions $u^* \in W_{loc}^{1,2}(\mathbb{R}^n)$ into finding $u \in W^{1,2}(\mathbb{R}/\Gamma')$, where $u^* = \alpha \cdot x + u$. Using Moser's idea, we establish periodic solutions corresponding to rotation vector $\alpha \in \mathbb{Q}^n$ by translating the effort of finding minimal and Birkhoff solutions into finding corresponding periodic configurations.

As in [15, 22], we look for minimal and Birkhoff solutions having the following form:

(5.9) $$u^*(\mathbf{i} + \mathbf{r}_i \mathbf{e}_i) = u^*(\mathbf{i}) + \mathbf{s}_i, \quad 1 \leq i \leq n,$$

where $\mathbf{r} \in \mathbb{N}^n$ and $\mathbf{s} \in \mathbb{Z}^n$. By (5.9) for $1 \leq i \leq n$ and $k \in \mathbb{Z}$,

(5.10) $$|u^*(\mathbf{i} + k\mathbf{r}_i \mathbf{e}_i) - \alpha \cdot (\mathbf{i} + k\mathbf{r}_i \mathbf{e}_i)| = |u^*(\mathbf{i}) + k\mathbf{s}_i - \alpha \cdot \mathbf{i} - \alpha_i k \mathbf{r}_i|.$$

Since $u^*$ is Birkhoff, by Lemma 2.5, there exists an $M > 0$ such that

$$|u^*(\mathbf{i}) - \alpha \cdot \mathbf{i}| \leq M$$

for all $\mathbf{i} \in \mathbb{Z}^n$. Noting $k$ is arbitrary in (5.10), we have

(5.11) $$\alpha_i = \mathbf{s}_i/\mathbf{r}_i, \quad 1 \leq i \leq n.$$

Thus setting $u = u^* - \alpha \cdot \mathbf{i}$ leads to

(5.12) $$u(\mathbf{i} + \mathbf{r}_i \mathbf{e}_i) = u^*(\mathbf{i} + \mathbf{r}_i \mathbf{e}_i) - \alpha \cdot (\mathbf{i} + \mathbf{r}_i \mathbf{e}_i) = u^*(\mathbf{i}) + \mathbf{s}_i - \alpha \cdot \mathbf{i} - \alpha_i \mathbf{r}_i = u(\mathbf{i}),$$

or $u \in \Gamma_0(\mathbf{r})$ ($\mathbf{r} = (\mathbf{r}_1, \cdots, \mathbf{r}_n)$). So if $\alpha \in \mathbb{Q}^n$ has the form of (5.11) with $\mathbf{r}_i, \mathbf{s}_i$ relatively prime, then searching for solutions with form as in (5.9) transfers to find periodic solutions $u$ defined as above.

For fixed $\alpha = (\alpha_1, \cdots, \alpha_n) \in \mathbb{Q}^n$, taking $\alpha_i = \mathbf{s}_i/\mathbf{r}_i$ with $\mathbf{s}, \mathbf{r}$ satisfying the above conditions (relatively prime, etc.), define

$$\begin{aligned}\Gamma_0^{\mathbf{r},\mathbf{s}} &:= \{u^* \in R^{\mathbb{Z}^n} \,|\, u^* \text{ satisfitis } (5.9)\}\\ &= \{u + \alpha \cdot \mathbf{i} \,|\, u \in \Gamma_0(\mathbf{r})\}\\ &= \Gamma_0(\mathbf{r}) + \alpha \cdot \mathbf{i},\end{aligned}$$

and

(5.13) $$c_0^{\mathbf{r},\mathbf{s}} := \inf_{u \in \Gamma_0^{\mathbf{r}}} J_0^{\mathbf{r}}(u + \alpha \cdot \mathbf{i}),$$

where $\Gamma_0(\mathbf{r})$, $J_0^{\mathbf{r}}$ are defined in Section 3.

Set
$$\mathcal{M}_0^{\mathbf{r},\mathbf{s}} := \{u + \alpha \cdot \mathbf{i} \,|\, u \in \Gamma_0^{\mathbf{r}} \text{ and } J_0^{\mathbf{r}}(u + \alpha \cdot \mathbf{i}) = c_0^{\mathbf{r},\mathbf{s}}\}.$$

Here $\mathcal{M}_0^{\mathbf{r},\mathbf{s}}$ corresponds to $\mathcal{M}_0$, which is defined for the case of $\alpha = \mathbf{0}$.

The following theorem is contained in [16]. It can also be proved as in [22, Theorem 5.27] and we omit the proof here.

**Theorem 5.7.** *(1) $\mathcal{M}_0^{\mathbf{r},\mathbf{s}} \neq \emptyset$ is an ordered set.*
*(2) Any $u^* = u + \alpha \cdot \mathbf{i} \in \mathcal{M}_0^{\mathbf{r},\mathbf{s}}$ is a minimal and Birkhoff solution of (2.1).*
*(3) For $\mathbf{k} \in \mathbb{N}^n$ and $\mathbf{t} \in \mathbb{Z}^n$, set $\hat{\mathbf{k}}(\mathbf{t}) = (\mathbf{k}_1\mathbf{t}_1, \cdots, \mathbf{k}_n\mathbf{t}_n)$. Then $M_0^{\hat{\mathbf{k}}(\mathbf{r}), \hat{\mathbf{k}}(\mathbf{s})} = \mathcal{M}_0^{\mathbf{r},\mathbf{s}}$ and*

$$c_0^{\hat{\mathbf{k}}(\mathbf{r}), \hat{\mathbf{k}}(\mathbf{s})} = (\prod_1^n \mathbf{k}_i) c_0^{\mathbf{r},\mathbf{s}}.$$

Now we consider $\mathcal{M}_1^\alpha$. To be brief we set $\Gamma_0^\alpha = \Gamma_0^{\mathbf{r},\mathbf{s}}$, $c_0^\alpha = c_0^{\mathbf{r},\mathbf{s}}$, and $\mathcal{M}_0^\alpha = \mathcal{M}_0^{\mathbf{r},\mathbf{s}}$. Assume that

$(*_0^\alpha)$ there are adjacent $v_0^\alpha, w_0^\alpha \in \mathcal{M}_0(\mathbf{r})$ with $v_0^\alpha < w_0^\alpha$.

The elements of $\mathcal{M}_1^\alpha$ have the form $U^* = U + \alpha \mathbf{i}$ with $U$ heteroclinic in $\mathbf{i}_1$ from $v_0^\alpha$ to $w_0^\alpha$. To use minimization argument, we introduce the suitable set and associated functional. Letting

$$\mathbb{T}_\alpha^{n-1} = \mathbb{Z}/\mathbf{r}_2 \times \cdots \times \mathbb{Z}/\mathbf{r}_n$$

and

$$\{i\mathbf{r}_1\} \times \{0, \cdots, \mathbf{r}_2 - 1\} \times \cdots \times \{0, \cdots, \mathbf{r}_n - 1\} =: \mathbb{T}_i^\alpha$$

replace $(\mathbb{Z}/\{1\})^{n-1}$ and $\mathbf{T}_i$, respectively. Then we can generalize $\hat{\Gamma}_1$, $J_{1,i}$, and $J_1$, etc. by

$$\hat{\Gamma}_1^\alpha := \hat{\Gamma}_1^\alpha(v, w) := \{u \in \mathbb{R}^{\mathbb{Z} \times \mathbb{T}_\alpha^{n-1}} \,|\, u \text{ lies between } v \text{ and } w\},$$

$$J_{1,i}^\alpha(u) := \sum_{\mathbf{j} \in \mathbb{T}_i^\alpha} S_{\mathbf{j}}(u + \alpha \cdot \mathbf{i}) - c_0^\alpha,$$

and
$$J_1^\alpha(u) := \liminf_{\substack{p\to -\infty \\ q\to\infty}} J_{1;p,q}^\alpha(u),$$

etc. Finally, define
$$\Gamma_1^\alpha := \Gamma_1^\alpha(v,w) := \{u \in \hat{\Gamma}_1^\alpha \mid \|u-v\|_{\mathbb{T}_i^\alpha} \to 0, i\to -\infty; \|u-w\|_{\mathbb{T}_i^\alpha} \to 0, i\to\infty\}.$$

Using these notations and the arguments of the previous sections, we obtain the results of the case $\alpha \in \mathbb{Q}^n$. For example, the next theorem is a new version of Theorem 3.13.

**Theorem 5.8.** *If $s \in C^2(\mathbb{R}^{B_0^r}, \mathbb{R})$ satisfies (S1)-(S3) and $(*_0^\alpha)$ holds, then there is a solution of the form $U_1^\alpha + \alpha \cdot \mathbf{i}$, where $U_1^\alpha \in \mathcal{M}_1^\alpha := \mathcal{M}_1^\alpha(v_0^\alpha, w_0^\alpha) := \{u \in \Gamma_1^\alpha(v_0^\alpha, w_0^\alpha) \mid J_1^\alpha(u) = c_1^\alpha\}$. Moreover, $\mathcal{M}_1^\alpha$ is an ordered set and the elements of $\mathcal{M}_1^\alpha + \alpha \cdot \mathbf{i}$ are solutions of (2.1), and any $U \in \mathcal{M}_1^\alpha$ is strictly $\mathbf{r}_1$-monotone in $\mathbf{i}_1$.*

**Acknowledgments** The authors would like to thank Xifeng Su for helpful discussions. W.-L. Li would like to thank Liang Jin for a lot of encouragement. X. Cui is supported by the National Natural Science Foundation of China (Grants 11571166, 11631006, 11790272), the Project Funded by the Priority Academic Program Development of Jiangsu Higher Education Institutions (PAPD) and the Fundamental Research Funds for the Central Universities.

Wen-Long Li, School of Mathematics, Sun Yat-Sen University, Guangzhou, 510275, P. R. China

*E-mail address*: liwenlongchn@gmail.com

Xiaojun Cui, Department of Mathematics, Nanjing University, Nanjing, 210093, P. R. China

*E-mail address*: xcui@nju.edu.cn